\documentclass[aop,seceqn,MSNbibl,dvips]{arximspdf}
\usepackage{mathbh}
\usepackage{graphicx}

\doi{10.1214/09-AOP520}
\volume{38}
\issue{4}
\pubyear{2010}
\firstpage{1690}
\lastpage{1716}

\makeatletter

\newproclaim{definition}{Definition}[section]
\newtheorem{proposition}[definition]{Proposition}
\newtheorem{theorem}[definition]{Theorem}
\newproclaim{remark}[definition]{Remark}
\newtheorem{lemma}[definition]{Lemma}
\newtheorem{corollary}[definition]{Corollary}

\makeatother

\begin{document}
\begin{frontmatter}

\title{Infinite rate mutually catalytic branching\protect\thanksref{T1}}
\runtitle{Infinite rate mutually catalytic branching}

\thankstext{T1}{Supported in part by the German Israeli Foundation
Grant G-807-227.6/2003.}

\begin{aug}
\author[A]{\fnms{Achim} \snm{Klenke}\corref{}\ead[label=e1]{math@aklenke.de}} and
\author[B]{\fnms{Leonid} \snm{Mytnik}\ead[label=e2]{leonid@ie.technion.ac.il}}
\runauthor{A. Klenke and L. Mytnik}
\affiliation{University of Mainz and Technion Haifa}
\address[A]{Institut f{\"u}r Mathematik\\
Johannes Gutenberg-Universit{\"a}t Mainz\\
Staudingerweg 9\\
55009 Mainz\\
Germany\\
\printead{e1}}
\address[B]{Faculty of Industrial Engineering\\
\quad and Management\\
Technion---Israel Institute of Technology\\
Haifa 32000\\
Israel\\
\printead{e2}}
\end{aug}

\received{\smonth{9} \syear{2008}}
\revised{\smonth{10} \syear{2009}}

%
\begin{abstract}
Consider the mutually catalytic branching process with finite branching
rate $\gamma$. We show that as $\gamma\to\infty$, this process
converges in finite-dimensional distributions (in time) to a certain
discontinuous process. We give descriptions of this process in terms of
its semigroup in terms of the infinitesimal generator and as the solution of a
martingale problem. We also give a strong construction in terms of a
planar Brownian motion from which we infer a path property of the
process.

This is the first paper in a series or three, wherein we also construct
an interacting version of this process and study its long-time
behavior.
\end{abstract}

%

\end{frontmatter}

\section{Introduction and main results}\label{A1}
\subsection{Motivation}
\label{A1.1}

In \cite{DawsonPerkins1998}, Dawson and Perkins introduced a
population dynamic model of two populations that live on a countable
site space $S$. The
individuals migrate between sites and, at any given site, perform a
critical branching process with a branching rate proportional to the local size of the population of
the respective other type.

More precisely, Dawson and Perkins considered the system of coupled
stochastic differential equations (SDEs) (taking nonnegative values)
%
%
\begin{equation}
\label{E1.1}
dY_{i,t}(k)=({\mathcal{A}}Y_{i,t})(k)\,dt + \sqrt{\gamma
Y_{1,t}(k)Y_{2,t}(k)}
\,dW_{i,t}(k),\qquad  i=1,2, k\in S.\hspace*{-33pt}
\end{equation}

Here, ${\mathcal{A}}(k,l)=a(k,l)-\mathbh{1}_{\{k\}}(l)$ is the
$q$-matrix of a Markov
chain on $S$ with symmetric jump kernel $a$, $(W_i(k), k\in S,
i=1,2)$ is an independent family of Brownian motions and $\gamma\geq
0$ is a parameter.

Dawson and Perkins showed that there exists a unique weak solution of
this SDE
taking values in a suitable subspace of $([0,\infty)^2)^S$ with some growth
condition. Furthermore, this process is a strong Markov process. While
existence of a weak solution is rather standard due to the procedure
proposed by
Shiga and Shimizu~\cite{ShigaShimizu1980}, weak
uniqueness was shown using a certain self-duality of the process established
in \cite{Mytnik1998}. We will describe the duality in detail below, in
(\ref{E2.4}).

A main result of Dawson and Perkins is a dichotomy in the long-time
behavior of the solutions depending on whether ${\mathcal{A}}$ is
recurrent or
transient (assuming some mild regularity condition on ${\mathcal
{A}}$). For
recurrent ${\mathcal{A}}$ (fulfilling the regularity assumption), the types
segregate, while for transient ${\mathcal{A}}$, there is coexistence
of types.
More precisely, let
\[
M_{i,t}=\sum_{k\in S}Y_{i,t}(k)
\]
denote the total mass processes ($i=1,2$) and assume that $M_{1,0},
M_{2,0}<\infty$. Then $M_1$ and $M_2$ are continuous orthogonal
nonnegative $L^2$-martingales. Let $M_{i,\infty}=\lim_{t\to\infty
}M_{i,t}$ denote the almost sure limit. Dawson and Perkins show that
$\mathbf{E}[M_{1,\infty} M_{2,\infty}]=0$ if ${\mathcal{A}}$ is
recurrent and
$\mathbf{E}[M_{1,\infty} M_{2,\infty}]=M_{1,0}M_{2,0}$ if ${\mathcal
{A}}$ is
transient. Furthermore, in the recurrent case, the joint distribution
of $(M_{1,\infty},M_{2,\infty})$ equals $Q_{(M_{1,0},M_{2,0})}$,
where, for $x\in[0,\infty)^2$, $Q_x$ is the harmonic measure of
planar Brownian motion in $[0,\infty)^2$. That is, if $B=(B_1,B_2)$ is
a Brownian motion in ${\mathbb{R}}^2$ started at $x$ and $\tau=\inf\{
t>0\dvtx
B_t\notin(0,\infty)^2\}$, then $Q_x$ is the probability measure on
\[
E:=[0,\infty)^2\setminus(0,\infty)^2
\]
given by
%
%
\begin{equation}
\label{e1.2}
Q_x=\mathbf{P}_x[B_\tau\in\bolds{\cdot}].
\end{equation}
The explicit form of the densities of $Q_x$ can be found in (\ref{E2.5}).

Via the self-duality of the mutually catalytic branching process, its total
mass behavior for finite initial conditions provides information on the local
behavior if the initial condition is infinite and sufficiently homogeneous.
For $x\in[0,\infty)^2$, let $\underline x$ denote the state in
$([0,\infty)^2)^S$ with $\underline x_i(k)=x_i$ for all $k\in S$, $i=1,2$.
Assume that $Y_0=\underline x$. Then
\[
\lim_{t\to\infty}\mathbf{P}_{\underline x}[Y_{1,t}(0)Y_{2,t}(0)>0]>0,
\]
if ${\mathcal{A}}$ is transient, that is, types can coexist locally.
On the
other hand,
for recurrent~${\mathcal{A}}$, the distribution of $Y_t$ converges
weakly to
$\int\delta_{\underline y}Q_x(dy)$,
that is, to a spatially homogeneous point $\underline y$, where $y$ is sampled
according to the distribution $Q_x$. Hence, in the recurrent case, the
two types
segregate locally and form clusters. The assumption that the initial
point is
constant can be weakened to an ergodic random initial condition (see
\cite{CoxKlenkePerkins2000}).

The starting point for this work was the wish to obtain a quantitative
description
of the cluster growth in the recurrent case. We will only briefly
describe the
heuristics. Dawson and Perkins also constructed a version of their
process in
continuous space ${\mathbb{R}}$ instead of $S$ as the solution of a
stochastic partial
differential equation
%
%
\begin{equation}
\label{e1.3}\qquad
\frac{dY_{i,t}(r)}{dt}=\Delta Y_{i,t}(r) + \sqrt{\gamma
Y_{1,t}(r)Y_{2,t}(r)} \dot{W}_i(t,r), \qquad r\in{\mathbb{R}}, i=1,2,
\end{equation}
where $\dot{W}_1$ and $\dot{W}_2$ are independent space--time white
noises and $\Delta$ is the Laplace operator. As $\Delta$ on ${\mathbb
{R}}$ is
recurrent, types also segregate here. Now, due to Brownian scaling, if
we denote by $Y^\gamma$ the solution of (\ref{e1.3}) with that given
value of $\gamma$, then we obtain
%
%
\begin{equation}
\label{E1.4}
\mathbf{P}_{\underline x} \bigl[\bigl(Y^\gamma_{T} \bigl(r\sqrt{T}\bigr)
\bigr)_{r\in{\mathbb{R}}
}\in\bolds{\cdot} \bigr]
=\mathbf{P}_{\underline x} [ (Y^{\gamma T}_1(r) )_{r\in
{\mathbb{R}}}\in
\bolds{\cdot} ].
\end{equation}
Equation (\ref{E1.4}) shows that clusters of $Y_{1,T}$ grow like
$\sqrt
{T}$ and
that a better understanding of the precise cluster formation can be
obtained by letting $\gamma\to\infty$ for fixed time. Hence, we aim
to construct a model $X$, that is, in some sense, the limit of $Y^\gamma
$ as $\gamma\to\infty$.

In this paper, we construct $X$ in the simple case where $S$ is a
singleton and where the migration between colonies is replaced by an
interaction with a time-invariant mean field. This is a first step
toward the investigation of the model involving infinitely many sites.
We give characterizations of the process $X$ via an
infinitesimal generator, as the solution of a well-posed martingale
problem and as
the limit of $Y^\gamma$ as $\gamma\to\infty$. Finally, we give a
strong construction of the process via a time-changed planar Brownian
motion. This will also serve to derive path properties.

In two forthcoming papers, we construct the infinite rate process on a
countable site space $S$ via a stochastic differential equation with
jump-type noise and give a characterization via a martingale problem
\cite{KM2}. Furthermore, we will investigate the long-time behaviour
and give conditions for segregation and for coexistence of types
\cite{KM3}. An alternative construction via a Trotter product approach is
carried out in \cite{KO} and \cite{Oeler2008}.

\subsection{Results}
\label{A1.2}

We now describe the one-colony process which is the subject of
investigation of this paper. Assume that $S$ is a singleton and that
immigration and emigration come from and go to some colony that is
thought to be infinitely big and whose effective population size (for
immigration) is $\theta\in[0,\infty)^2$.
Furthermore, let $c\geq0$ be the rate of migration.
Hence, we consider the solution $Y=Y^{\gamma,c,\theta}$ of the
stochastic differential equation
%
%
\begin{equation}
\label{E1.5}
dY_{i,t}=c (\theta_i-Y_{i,t} ) \,dt + \sqrt{\gamma
Y_{1,t}Y_{2,t}} \,dW_{i,t}, \qquad i=1,2.
\end{equation}
This model can be thought of as a version of the model defined in (\ref
{E1.1}) where the migration between colonies is replaced by an
interaction with a time-invariant mean field $\theta$ or with an
infinitely large reservoir whose types have proportions $\theta_1$ and
$\theta_2$. (In fact, in \cite{CoxDawsonGreven2004} it was shown
(Proposition 1.1) that $Y^{\gamma,c,\theta}$ arises as the
McKean--Vlasov limit of solutions of (\ref{E1.1}) with symmetric
interaction on a complete graph $S$.) More formally, the interaction
term ${\mathcal{A}}Y$ is replaced by a drift $c (\theta
_i-Y_{i,t} )$.
It is this simplification of the interaction that allows for a
tractable exposition in this article.
Note that as $t\to\infty$, the process without drift ($c=0$)
converges almost surely to some random $x\in E$. Hence, in the case
$c=0$, if we let $\gamma\rightarrow\infty$, then the limiting
process would be trivial: if it starts at $x\in E$, then it stays at
$x$ forever.
See Section \ref{A2} for a more detailed description of
the process $Y$ solving (\ref{E1.5}) (finite $\gamma$ process).

On a heuristic level, as the stochastic term in (\ref{E1.5}) defines an
isotropic two-dimensional diffusion, that is, a time-transformed planar
Brownian motion, if we let $\gamma\to\infty$, then we should end up
with a process where the stochastic part is a planar Brownian motion at
infinite speed, stopped when it reaches the boundary of the upper-right
quadrant. That is, the limiting process $X$ should be a Markov process
with values in $E$. When $x$ is the current state and the drift moves
it to $x+c(\theta-x)\,dt$, this point should instantaneously be replaced
by a random point chosen according to $Q_{x+c(\theta-x)\,dt}$. We will,
in fact, be able to describe this infinitesimal dynamics both in terms
of a martingale problem and in terms of a generator of Markov
transition kernels. However, we first define $X$ via an explicit
transition semigroup and show that it is the limit of $Y^{\gamma
,c,\theta}$ as $\gamma\to\infty$.
Let
%
%
\begin{equation}
\label{E1.6}
C_l(E):= \Bigl\{f\dvtx E\to{\mathbb{C}}\mbox{ is cont. and }\lim_{u\to
\infty
}f(u,0)=\lim_{v\to\infty}f(0,v)\mbox{ is finite} \Bigr\}\hspace*{-28pt}
\end{equation}
equipped with the supremum norm $\|f\|_\infty={\sup_{x\in E}}|f(x)|$.
\begin{definition}
\label{D1.1}
Let $c\geq0$ and $\theta\in[0,\infty)^2$.
For $t\geq0$ and $x\in E$, define the stochastic kernel $p_t$ by
\[
p_t(x,\bolds{\cdot}):=p_t^{c,\theta}(x,\bolds{\cdot
}):=Q_{e^{-ct}x+(1-e^{-ct})\theta}.
\]
Define the contraction semigroup $\mathcal{S}=(\mathcal{S}_t)_{t\geq
0}$ on $C_l(E)$ by
\[
\mathcal{S}_tf(x)=\int_Ef(y) p_t(x,dy).
\]
The Markov process $X=X^{c,\theta}$ with state space $E$, c{\`a}dl{\`
a}g paths and transition kernels $(p_t)_{t\geq0}$ is called the
infinite rate mutually catalytic branching process (IMUB) with
parameters $(c,\theta)$.
\end{definition}

In order for this definition to make
sense, in Proposition \ref{P3.2}, we will show that $(\mathcal
{S}_t)_{t\geq
0}$ is, in fact, a Markov semigroup.
\begin{proposition}
\label{P1.2}
$X^{c,\theta}$ is a Feller process and has the strong Markov property.
It is ergodic and the unique invariant measure is $Q_\theta$.
\end{proposition}
\begin{pf}
The map $x\mapsto Q_x$ is continuous, hence $x\mapsto p_t(x,\bolds
{\cdot})$ is
also continuous, that is, $X^{c,\theta}$ is a Feller process. Since
$Q_x=\delta_x$ for $x\in E$, the semigroup $\mathcal{S}$ is strongly
continuous. Hence, by the general theory of Markov processes, there
exists a c{\`a}dl{\`a}g version of $X$ that is strong Markov (see,
e.g., \cite{RogersWilliams1994}, Chapters III.7 and 8).

Ergodicity and the explicit form of the invariant measure are trivial.
\end{pf}
\begin{theorem}[($X^{c,\theta}$ as an infinite rate process)]
\label{T1.3}
Assume that $Y^{\gamma,c,\theta}_0=X^{c,\theta}_0=x\in E$ for all
$\gamma\geq0$. As $\gamma\to\infty$, the finite-dimensional
distributions of $Y^{\gamma,c,\theta}$ converge to those of
$X^{c,\theta}$.
\end{theorem}

Note that in Theorem \ref{T1.3}, trivially, we do not have convergence
in the Skorohod path space, since continuous processes do not converge
to discontinuous processes in that topology.

In addition to the convergence of the finite-dimensional distributions,
we also have convergence of the $p$th moments for $p\in[1,2)$ [but not
for $p=2$, of course, since for $x\in(0,\infty)^2$, the measure $Q_x$
does not possess finite second moments, as can be easily derived from
its density formula (\ref{E2.5})]. Hence, on a suitable probability
space, we have $L^p$-convergence of $Y^{\gamma,c,\theta}$ to
$X^{c,\theta}$.
\begin{theorem}[($L^p$-convergence)]
\label{T1.4}
Assume that $Y^{\gamma,c,\theta}_0=X^{c,\theta}_0=x\in E$ for all
$\gamma\geq0$ and let $p\in[1,2)$, $t\geq0$.

\begin{longlist}
\item
For every $\gamma\geq0$ and $i=1,2$, we have
\[
\mathbf{E}_x [(Y^{\gamma,c,\theta}_{i,t})^p ]\leq
\mathbf{E}_x [(X^{c,\theta}_{i,t})^p ]<\infty.
\]

\item
On a suitable probability space, for $i=1,2$, we have
\[
Y^{\gamma,c,\theta}_{i,t}\stackrel{\gamma\to\infty
}{\longrightarrow}X^{c,\theta
}_{i,t} \qquad\mbox{in }L^p.
\]
\end{longlist}
\end{theorem}

It can be seen from the proofs of Theorems \ref{T1.3} and \ref{T1.4}
that the statements of these theorems also hold for $Y^{\gamma
,c,\theta}_0=x\in[0,\infty)^2$ and $t>0$ if we replace $X^{c,\theta
}_0$ by a random point chosen according to $Q_{x}$.
\begin{remark}[(Trotter product approach)]
\label{R3.3}
While
in the one-colony case considered in this paper, it is easy to
explicitly write down the semigroup for the infinite
rate mutually catalytic branching process $X^{c,\theta}$, it is less obvious
how to construct an interacting version of the process on a countable site
space. One possibility is the Trotter product approach that is used in
\cite{KO} and \cite{Oeler2008}. Here, we briefly sketch it for
$X^{c,\theta}$.

In the classical setting, the Trotter product approach works as follows.
In order to construct a solution $Y^{\gamma,c,\theta}$ of (\ref
{E1.5}), in time intervals of length ${\varepsilon}$, one could
alternate between
a solution of the pure drift equation ($\gamma=0$) and the pure
stochastic noise equation ($c=0$). As ${\varepsilon}\downarrow0$,
this process
converges to a solution of (\ref{E1.5}).

If we let $\gamma\to\infty$, then the noise term results in an
instantaneous jump to
a point in $E$ chosen according to $Q_y$, where $y$ is the value of $Y$
at the
end of the preceding ``drift interval.'' More formally, let $(\xi
(k,x), k\in{\mathbb{N}},
x\in[0,\infty)^2)$ be an independent family of $E$-valued random variables
with distribution ${\mathcal{L}}[\xi(k,x)]=Q_x$. For $t\in
[k{\varepsilon},(k+1){\varepsilon})$,
let $X^{\varepsilon}_t$
be the solution of the differential equation
\[
dX^{\varepsilon}_t=c(\theta-X_t)\,dt,
\]
that is,
\[
X^{\varepsilon}_t=e^{-c(t-k{\varepsilon})}X^{\varepsilon
}_{k{\varepsilon}}+ \bigl(1-e^{-c(t-k{\varepsilon})}
\bigr)\theta.
\]
Let
\[
X^{{\varepsilon}}_{(k+1){\varepsilon}-}:=\lim_{t\uparrow
(k+1){\varepsilon}}X^{\varepsilon}_t=e^{-c{\varepsilon}
}X^{\varepsilon}_{k{\varepsilon}}+(1-e^{-c{\varepsilon}})\theta
\]
and define
\[
X^{\varepsilon}_{(k+1){\varepsilon}}=\xi \bigl(k+1,X^{\varepsilon
}_{(k+1){\varepsilon}-} \bigr).
\]
One can prove that
$X^{\varepsilon}$ converges in distribution in the Skorohod topology
on the
space of
c{\`a}dl{\`a}g paths to $X^{c,\theta}$ (see \cite{KO} and
\cite{Oeler2008}).
\end{remark}

While, in Definition \ref{D1.1}, we gave an explicit formula for the
transition kernels of $X$, it is also interesting to characterize the
process $X$ via its infinitesimal dynamics. In Section \ref{A5}, we
investigate the generator $\bar{\mathcal{G}}$ of the semigroup
$\mathcal{S}$. For a
certain class $C_l^2(E)\subset C_l(E)$ of smooth functions $f$ (see
Definition \ref{D5.1}), we give an explicit formula for $\bar
{\mathcal{G}}f$
as an integro-differential operator. Using the classical Hille--Yoshida
theorem, we show that the restricted operator ${\mathcal{G}}=\bar
{\mathcal{G}}|_{C_l^2(E)}$
uniquely defines $(\mathcal{S}_t)_{t\geq0}$ (Theorem \ref
{T5.3}). Furthermore, we show that ${\mathcal{G}}$ restricted to an
even smaller
space $V$ of functions that
appear in the duality for $X$ still uniquely defines the process $X$
via a martingale problem (Theorem \ref{T5.4}).
To define $ {\mathcal{G}}$, it is crucial to study (for suitable
functions $f$)
the limit
\[
\lim_{t\downarrow0}t^{-1}\bigl(\mathcal{S}_tf(x)-f(x)\bigr)=\lim
_{{\varepsilon}\to0}{\varepsilon}
^{-1} \biggl(\int f \,dQ_{x+\varepsilon c(\theta-x)}-f(x) \biggr),
\]
which will also clarify the jump structure of the process $X$.
The description of the exact form of the operator ${\mathcal{G}}$ and the
precise statements of the theorems are a bit
technical, so these are deferred to Section \ref{A5}.

While, for Proposition \ref{P1.2}, we used general construction
principles of Markov processes, here, we provide an explicit strong
construction of the process $X$
in terms of a given planar Brownian motion $B$. This construction also
allows certain path properties to be investigated.

Assume $B_0=0$. For $z\in{\mathbb{R}}^2$, we write
\[
[z,\infty)=[z_1,\infty)\times[z_2,\infty)
\]
for the rectangular cone northeast of $z$.
For $x\in[0,\infty)^2$, let
%
%
\begin{equation}
\label{E1.7}
\tau_x:=\inf \{t>0\dvtx B_t\notin[-x,\infty) \}
\end{equation}
and
%
%
\begin{equation}
\label{E1.8}
D_x:=B_{\tau_x}+x \in E.
\end{equation}
For $x,y\in{\mathbb{R}}^2$, we write $y\leq x$ if $x\in[y,\infty)$,
that is,
if $y_1\leq x_1$ and $y_2\leq x_2$. For $x\in[0,\infty)^2$, we define
the $\sigma$-algebra
%
%
\begin{equation}
\label{E1.8b}
{\mathcal{F}}^D_x=\sigma(D_y\dvtx y\leq x).
\end{equation}
In Lemma \ref{L3.1}, we will show that $D$ is a Markov process with
respect to $({\mathcal{F}}^D_x)_{x\in[0,\infty) ^2}$.

Let $\bar\theta\dvtx[0,\infty)\to[0,\infty)^2$ and $\bar c\dvtx[0,\infty
)\to[0,\infty)$ be measurable and locally integrable. For $0\leq
s\leq t$, define
%
%
\begin{equation}
\label{E1.9}
C(s,t)=\exp \biggl(-\int_s^t \bar c(r) \,dr \biggr) \quad\mbox{and}\quad
\Xi(s,t)=\int_s^t\frac{\bar\theta(r)}{C(0,r)} \,dr.
\end{equation}
\begin{theorem}
\label{T1.5}
Let $x\in E$ and define the process $X^{\bar c,\bar\theta}$ by
\[
X^{\bar c,\bar\theta}_t=C(0,t)D_{x+\Xi(0,t)}, \qquad t\geq0.
\]
Then $X^{\bar c,\bar\theta}$ is a time-inhomogeneous Markov process
on $E$ with c{\`a}dl{\`a}g paths and with transition probabilities
%
%
\begin{equation}
\label{E1.10}
p_{s,t}(z,\bolds{\cdot})=Q_{C(s,t)z+C(0,t)\Xi(s,t)}
\qquad\mbox{for }0\leq s<t, z\in E.
\end{equation}
In particular, for $\bar\theta\equiv\theta\in[0,\infty)^2$ and
$\bar c\equiv c>0$,
%
%
\begin{equation}
\label{E1.11}
X^{c,\theta}_t=e^{-ct}D_{x+(e^{ct}-1)\theta}
\end{equation}
is an infinite rate mutually catalytic branching process
%
%
\begin{figure}[b]

\includegraphics{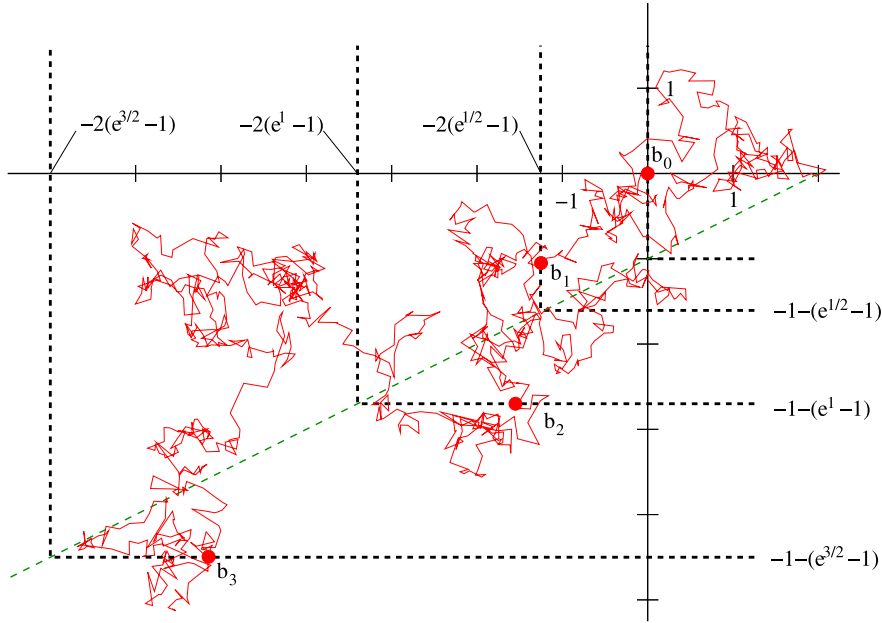}

\caption{Strong construction of $X^{1/2,(2,1)}$ with
$X_0=x=(0,1)$ via a planar Brownian motion. Here
$X^{1/2,(2,1)}_t=e^{-t/2}((0,1)+b_t+(2,1)(e^{t/2}-1))$ for
$t=0, 1, 2, 3$.}
\label{figure1}
\end{figure}
with parameter $(c,\theta)$, see Figure \ref{figure1}.
\end{theorem}

It is tempting to use this strong construction of $X^{\bar c,\bar
\theta}$ in order to define an interacting
version of the infinite rate mutually catalytic branching process
on a countable site space $S$, where $c\theta_k(t)$ at site $k\in S$
reflects the
migration from neighboring sites to $k$. However, in this paper, we do
not pursue this topic. Rather, we use the strong construction in order
to derive a path property of $X^{c,\theta}$ via a result of Le Gall
and Meyre \cite{LeGallMeyre1992} on the cone points of planar Brownian motion.

Recall that a measurable set $A\subset E$ is called \textit{polar} for
$X^{c,\theta}$ if for all $x\in E$, we have
\[
\mathbf{P}_x[X^{c,\theta}_t\in A\mbox{ for some }t>0]=0.
\]
\begin{theorem}
\label{T1.6}
The point $0\in E$ is polar for $X^{c,\theta}$.
\end{theorem}

\subsection{Organization of the paper}
\label{A1.3}

In Section \ref{A2}, we give a detailed description of the duality for
the process with finite branching rate. In Section \ref{A3}, we
establish a similar duality for the infinite rate process and use it in
order to show the convergence in Theorems \ref{T1.3} and \ref{T1.4}.
In Section \ref{A4}, we justify the strong construction of
Theorem \ref{T1.5} and also prove Theorem \ref{T1.6}. Finally, in
Section \ref{A5}, we describe the infinite rate process in terms of
its infinitesimal dynamics and state and prove the theorem on the
construction via the Hille--Yoshida theory (Theorem \ref{T5.3}) and
via a martingale problem (Theorem \ref{T5.4}).

\section{Duality of the finite $\gamma$ process}
\label{A2}

A major tool for the investigation of mutually catalytic branching
processes is a self-duality for the process. As it turns out to be
crucial also for the limiting case of infinite branching rate ($\gamma
=\infty$), we describe this duality here in more detail.
For $x=(x_1,x_2)$ and $y=(y_1,y_2)\in{\mathbb{R}}^2$, we introduce
the \textit{lozenge product}
%
%
\begin{equation}
\label{E2.1}
x \diamond y := -(x_1+x_2)(y_1+y_2) + i(x_1-x_2)(y_1-y_2)
\end{equation}
(with $i=\sqrt{-1}$) and define
%
%
\begin{equation}
\label{E2.2}
F(x, y)=\exp(x\diamond y).
\end{equation}
Note that $x\diamond y=y\diamond x$. Furthermore, define the ``scalar
product''
%
%
\begin{equation}
\label{E2.3}
\langle x,y\rangle=x_1y_1+x_2y_2 \qquad\mbox{for }x,y\in[0,\infty)^2.
\end{equation}

For $x=(x(k))_{k\in S}$ and $y=(y(k))_{k\in S}$, we write
\[
H(x, y)=\exp \biggl(\sum_{k\in S} x(k)\diamond y(k) \biggr).
\]
If $Y$ is the process defined in (\ref{E1.1}) started in state $y$ and
$\tilde Y$ is the process started in some suitable $\tilde y$ (such
that all sums are finite), then the duality reads (see \cite{Mytnik1998}, equation
(2.5))
%
%
\begin{equation}
\label{E2.4}
\mathbf{E}_y [H (Y_t,\tilde y ) ]=\mathbf{E}_{\tilde
y} [H
(y,\tilde Y_t ) ].
\end{equation}
In fact, this duality also holds for asymmetric ${\mathcal{A}}$ if
$\tilde Y$ is a
solution of (\ref{E1.1}) with ${\mathcal{A}}$ replaced by its
transpose ${\mathcal{A}}^*$.
As this mixed Laplace and Fourier transform $H$ is measure determining
(\cite{Mytnik1998}, Lemma 2.5), the
duality yields uniqueness of the solutions of (\ref{E1.1}).
Furthermore, it
provides a tool for translating local properties of the solutions into global
properties and vice versa.
If $x=(u,v)\in(0,\infty)^2$, then the harmonic measure $Q_x$ [recall
(\ref{e1.2})] has a one-dimensional Lebesgue density on
\[
E:= \bigl([0,\infty)\times\{0\} \bigr)\cup \bigl(\{0\}\times[0,\infty
) \bigr)
\]
that can be computed explicitly
%
%
\begin{equation}
\label{E2.5}
Q_{(u,v)} (d(\bar u,\bar v) )=
\cases{
\displaystyle\frac4\pi \frac{uv \bar u}{
4u^2v^2+ (\bar u^2+ v^2-u^2 )^2} \,d\bar u, &\quad  if
$\bar v=0$,\cr
\displaystyle\frac4\pi \frac{uv \bar v}{4u^2v^2+ (\bar
v^2+u^2-v^2 )^2} \,d\bar v, &\quad  if $\bar u=0$.}
\end{equation}
Furthermore, trivially we have
%
%
\begin{equation}
\label{E2.6}
Q_x=\delta_x \qquad\mbox{if }x\in E.
\end{equation}

We now turn to the situation of only one colony.
We consider the solution $Z=(Z_1,Z_2)$ of
%
%
\begin{equation}
\label{E2.7}
dZ_{i,t}=\sqrt{\gamma Z_{1,t}Z_{2,t}} \,dW_{i,t},\qquad  i=1,2,\qquad
Z_0=z\in[0,\infty)^2.
\end{equation}
By Theorem 1 of \cite{DawsonFleischmannXiong2005}, there is the unique
strong solution to the above equation.

Clearly, $Z_1$ and $Z_2$ are orthogonal $L^2$-martingales and hence
they converge almost surely to some random variable $Z_\infty
=(Z_{1,\infty},Z_{2,\infty})$. As $Z$ is an isotropic diffusion on
$[0,\infty)^2$, it is a time-transformed Brownian motion. Thus
$Z_\infty$ has the same distribution as a planar Brownian motion $B$
started at $z$ and stopped (at time $\tau$) upon leaving $(0,\infty
)^2$, that is [see (\ref{E2.5})],
\[
{\mathcal{L}}_z[Z_\infty]={\mathcal{L}}_z[B_\tau]=Q_z.
\]
(We denote by ${\mathcal{L}}_x[X_t]=\mathbf{P}_x[X_t\in\bolds{\cdot
}]=\mathbf{P}[X_t\in\bolds{\cdot}
\mid X_0=x]$ the distribution of the process $X$ at
time $t$ when started at $x$.) It is easy to see that, in fact,
\[
\tau^Z:=\inf\{t>0\dvtx Z_t\in E\}<\infty \qquad\mbox{almost surely},
\]
and that
\[
Z_t=Z_{\tau^Z} \qquad\mbox{for all }t>\tau^Z.
\]

Clearly, increasing $\gamma$ amounts to speeding up the process.
Hence, in the limit, we would have a process that instantaneously jumps
from $z$ to a random point (picked according to $Q_z$) and then stays
there. In order to obtain a more interesting limiting process, and with
a view toward interacting colonies, we introduce a drift term and
consider the equation (which was analyzed in more detail in \cite
{CoxDawsonGreven2004}, Propositions 1.1 and 1.2)
%
%
\begin{equation}
\label{E2.8}
dY_{i,t}=c (\theta_i-Y_{i,t} ) \,dt + \sqrt{\gamma
Y_{1,t}Y_{2,t}} \,dW_{i,t}, \qquad i=1,2.
\end{equation}
Here, $c\geq0$ and $\theta\in[0,\infty)^2$ are parameters of the process.
It is standard to show that (\ref{E2.8}) has a weak solution. Weak
uniqueness can be obtained via duality. We first outline the general
picture for the duality that comes from the interacting colonies case
and then give an explicit computation for our special situation.

Let us consider a two-colonies model with site space $S=\{1,2\}$, where
$Y$ is the size of the population at site 1 and the size of the
population at site 2 is constant and equals $\theta$. This amounts to
a migration matrix
%
%
\begin{equation}
\label{E2.9}
{\mathcal{A}}=\pmatrix{-c & c\cr
0&0}
\end{equation}
and to branching rates $\gamma(1)=\gamma$ (at site 1) and $\gamma
(2)=0$ (at site 2). Note that the approach of Dawson and Perkins does
not require that the branching rate be constant; neither does it
require that the migration matrix be symmetric or a $q$-matrix. (At
least if $S$ is finite---otherwise, certain regularity conditions have
to be imposed.) Dawson and Perkins use a duality with respect to a
process $\tilde Y$ with migration matrix ${\mathcal{A}}^*$ (the
transpose of
${\mathcal{A}}$) and with the same branching rates as $Y$ to show weak
uniqueness of $Y$.

Let us now construct the dual process explicitly. We will later use
this approach in order to construct a dual for the $\gamma=\infty$
limiting process. Let $\tilde y=(\tilde y(1),\tilde y(2))\in([0,\infty
)^2)^2$ and let $Z$ be the unique strong (by Theorem 1 of \cite
{DawsonFleischmannXiong2005}) $[0,\infty)^2$-valued solution of
%
%
\begin{equation}
\label{E2.10}
dZ_{i,t}=\sqrt{\gamma Z_{1,t}Z_{2,t}} \,dW_{i,t},\qquad  i=1,2,\qquad
Z_0=\tilde y(1).
\end{equation}
Define a process $\tilde Y$ on $([0,\infty)^2)^2$ by
%
%
\begin{equation}
\label{E2.11}
\tilde Y_t(1)=e^{-ct}Z_t \quad\mbox{and}\quad \tilde Y_t(2)=\tilde y(2)+
\int_0^tc e^{-cr}Z_r \,dr.
\end{equation}

Note that this $\tilde Y$ is a solution of (\ref{E1.1}) with $S=\{1,2\}
$, with site-dependent branching rate $\gamma(1)=\gamma$, $\gamma
(2)=0$ and with ${\mathcal{A}}$ from (\ref{E2.9}) replaced by
${\mathcal{A}}^*$.
In particular, $\tilde Y$ is a time-homogeneous Markov process. We also
get the time-homogeneous Markov property via an explicit computation:
\begin{eqnarray*}
\tilde Y_{t+s}
&=& \biggl(e^{-c(t+s)}Z_{t+s}, \tilde y(2)+\int_0^{t+s}c e^{-cr}Z_r
\,dr \biggr)\\
&=& \biggl(e^{-cs} (e^{-ct}Z_{t+s}), \tilde y(2)+\int_0^{t}c
e^{-cr}Z_r \,dr + \int_0^{s}c e^{-cr}(e^{-ct}Z_{t+r}) \,dr \biggr)\\
&=& \biggl(e^{-cs} Z'_s, \tilde y'(2)+ \int_0^{s}c e^{-cr}Z'_r
\,dr \biggr),
\end{eqnarray*}
where $Z'_r=e^{-ct}Z_{t+r}$ and $\tilde y'(2)=\tilde Y_t(2)=\tilde
y(2)+\int_0^{t}c e^{-cr}Z_r \,dr$. Clearly, $Z'$ has the distribution
of a solution of (\ref{E2.7}) with $\tilde y'(1):=Z'_0=\tilde Y_t(1)$.

For $x,x',y,y'\in[0,\infty)^2$, recall that
%
%
\begin{equation}
\label{E2.12}
H ((x,x'),(y,y') )=F(x,y) F(x',y').
\end{equation}
\begin{proposition}[(Duality)] Let $Y$ and $\tilde Y$ be defined by (\ref
{E2.8}) and (\ref{E2.11}), respectively. Then, for all $y\in[0,\infty
)^2$, $\tilde y\in([0,\infty)^2)^2$ and $t\geq0$, we have
\label{P2.1}
%
%
\begin{equation}
\label{E2.13}
\mathbf{E}_y [H ((Y_t,\theta),\tilde y ) ]
=\mathbf{E}_{\tilde y} [H ((y,\theta),\tilde Y_t ) ].
\end{equation}
In particular, if $Z$ is a solution of (\ref{E2.10}) with $Z_0=z\in
[0,\infty)^2$, then
%
%
\begin{equation}
\label{E2.14}
\mathbf{E}_y[F(Y_t,z)]
=\mathbf{E}_{z} \biggl[F(y,e^{-ct}Z_t) F \biggl(\theta,\int_0^tc
e^{-cr}Z_r
\,dr \biggr) \biggr].
\end{equation}
\end{proposition}

A similar duality was derived in \cite{CoxDawsonGreven2004}, Lemma 4.2.
Before we prove the proposition, we have to collect some properties of
the derivatives of $F$. We omit the proof of the following lemma.
\begin{lemma}[(Derivatives of the duality function)]
\label{L2.2}
Denote the partial derivatives of $F$ 
by
\[
\nabla_1 F(x,y):=\frac{d}{dx}F(x,y),\qquad
\nabla_2 F(x,y):=\frac{d}{dy}F(x,y)
\]
and define the Laplace operators $\Delta_1$ and $\Delta_2$ by
\[
\Delta_1 F(x,y):= \biggl[\frac{\partial^2}{\partial x_1^2}+\frac
{\partial^2}{\partial x_2^2} \biggr]F(x,y),\qquad
\Delta_2 F(x,y):= \biggl[\frac{\partial^2}{\partial y_1^2}+\frac
{\partial^2}{\partial y_2^2} \biggr]F(x,y).
\]
Then, for all $x,y,z\in[0,\infty)^2$, we have [recall (\ref{E2.1}) and
(\ref{E2.3})]
\begin{eqnarray*}
 \langle z, \nabla_1 F(x,y) \rangle &=& (z\diamond y)
F(x,y),\\
 \langle z, \nabla_2 F(x,y) \rangle &=& (z\diamond x)
F(x,y),\\
\Delta_1F(x,y)&=&8y_1y_2 F(x,y),\\
\Delta_2F(x,y)&=&8x_1x_2 F(x,y).
\end{eqnarray*}
\end{lemma}
\begin{pf*}{Proof of Proposition \ref{P2.1}}
We use It{\^o}'s formula and Lemma \ref{L2.2} to compute the
derivatives of both sides of (\ref{E2.13})
at $t=0$:
%
%
\begin{eqnarray}
\label{E2.15}
&&\frac{d}{dt}\mathbf{E}_y [H ((Y_t,\theta),\tilde y
) ]
|_{t=0}\nonumber\\
&&\qquad= \langle c(\theta-y), \nabla_1 F(y,\tilde
y(1)) \rangle F(\theta,\tilde y(2)) \nonumber\\[-8pt]\\[-8pt]
&&\qquad\quad{} +
\frac12\gamma  y_1y_2 \Delta_1 F(y,\tilde y(1)) F(\theta,\tilde
y(2))\nonumber\\
&&\qquad=H ((y,\theta),\tilde y )  [c(\theta-y)\diamond\tilde y(1)
 + 4\gamma y_1y_2\tilde y_1(1)\tilde y_2(1) ]\nonumber
\end{eqnarray}
and
%
%
\begin{eqnarray}
\label{E2.16}
&&\frac{d}{dt}
\mathbf{E}_{\tilde y} [H ((y,\theta),\tilde Y_t )
]|_{t=0}\nonumber\\
&&\qquad=F(\theta,\tilde y(2))  \biggl( \langle-c \tilde
y(1),
\nabla_2F(y,\tilde y(1)) \rangle\nonumber\\
&&\qquad\quad\hspace*{53.12pt}{}
+\frac{\gamma}{2}\tilde y_1(1)\tilde y_2(1)\Delta_2F(y,\tilde
y(1)) \biggr)\\
&&\qquad\quad{} + F(y,\tilde y(1))  \langle c\tilde y(1),\nabla_2 F(\theta
,\tilde y(2)) \rangle\nonumber\\
&&\qquad=H ((y,\theta),\tilde y )  [c(\theta-y)\diamond\tilde y(1)
 + 4\gamma y_1y_2\tilde y_1(1)\tilde y_2(1) ].\nonumber
\end{eqnarray}
Since the two derivatives coincide, (\ref{E2.13}) holds (see
Corollary 4.4.13 of \cite{EthierKurtz1986} with
$\alpha=\beta=0$). Equation (\ref{E2.14}) is a direct consequence of
(\ref{E2.13}).
\end{pf*}
\begin{corollary}
\label{C2.3}
Recall $Z$ defined by (\ref{E2.10}).

\begin{longlist}
\item Taking $c=0$, Proposition \ref{P2.1} implies that $Z$ is
self-dual:
\[
\mathbf{E}_x [F (Z_t,y ) ]=\mathbf{E}_y[F
(x,Z_t ) ] \qquad\mbox{for all }
x,y\in[0,\infty)^2, t\geq0.
\]

\item Letting $t\to\infty$ in \textup{(i)} and recalling that
${\mathcal{L}}
_x [Z_t ]\stackrel{t \rightarrow\infty}{\longrightarrow
}Q_x$, we get, by dominated convergence, the
duality relation for the harmonic measure:
\[
\int_E F(z,y)  Q_x(dz)= \int_E F(x,z) Q_y(dz) \qquad\mbox{for all
}x,y\in[0,\infty)^2.
\]

\item In particular (since $Q_x=\delta_x$ for $x\in E$ and
due to the symmetry of $F$), for all $x\in E$ and $y\in[0,\infty)^2$,
we have
\[
\int_E F(x,z) Q_y(dz) =F(x,y)
=F(y,x)=
\int_E F(z,x) Q_y(dz).
\]
\end{longlist}
\end{corollary}
\begin{corollary}
\label{C2.4}
\textup{(i)}
The family of functions
\[
{\mathcal{F}}_0= \{[0,\infty)^2\to{\mathbb{C}}\dvtx x\mapsto
F(x,y), y\in[0,\infty
)^2 \}
\]
is measure determining for $[0,\infty)^2$.
\smallskipamount=0pt
\begin{longlist}[(ii)]
\item[(ii)]
The vector space
%
%
\begin{equation}
\label{E2.17}
V:= \Biggl\{\sum_{m=1}^n\lambda_mF(\bolds{\cdot},z_m)\dvtx n\in
{\mathbb{N}},  \lambda
_1,\ldots,\lambda_n\in{\mathbb{C}}, z_1,\ldots,z_n\in E \Biggr\}
\end{equation}
spanned by ${\mathcal{F}}:=\{E\to{\mathbb{C}}\dvtx x\mapsto F(x,z),
z\in E\}$ is dense in $C_l(E)$.
In particular, ${\mathcal{F}}$
is measure determining for $E$.
\end{longlist}
\end{corollary}
\begin{pf}
Let $\mathcal{D}_0$ be the algebra generated by ${\mathcal{F}}_0$.
Clearly, ${\mathcal{F}}_0$
separates points of $[0,\infty)^2$, contains $1=F(\bolds{\cdot},0)$
and is
closed under multiplication and under complex conjugation since
$\overline{F(x,(y_1,y_2))}=F(x,(y_2,y_1))$. Hence, by the
Stone--Weierstrass theorem, $\mathcal{D}_0$ is dense in the space
$C_l([0,\infty)^2)$ of functions $[0,\infty)^2\to{\mathbb{C}}$ that are
continuous and have a limit at infinity. As ${\mathcal{F}}_0$ is
closed under
multiplication, $\mathcal{D}_0$ is the vector space spanned by
${\mathcal{F}}_0$ and
thus ${\mathcal{F}}_0$ is measure determining on $[0,\infty)^2$.

Let ${\mathcal{F}}_E=\{f|_{E}\dvtx f\in{\mathcal{F}}_0\}\supset
{\mathcal{F}}$ and let $\mathcal{D}
_E=\{f|_{E}\dvtx f\in\mathcal{D}_0\}$ denote the algebra generated by
${\mathcal{F}}
_E$. By the above argument, $\mathcal{D}_E\subset C_l(E)$ is dense.
Now, by
Corollary \ref{C2.3}(iii), an element $F(\bolds{\cdot},y)\in
{\mathcal{F}}_E$ can be
written as the integral $F(x,y)=\int F(x,z)Q_y(dz)$, where the
integrand functions are in ${\mathcal{F}}$. The integral can be approximated
(uniformly in $x$) by finite sums, that is, by elements of $V$. Hence,
$V$ is dense in $\mathcal{D}_E$ and thus also in $C_l(E)$.
\end{pf}

Apparently, $Y$ is ergodic and has a unique invariant distribution with
a Lebesgue density on $(0,\infty)^2$. Unlike for the analogous
one-dimensional equation
\[
dU_t=c(b-U_t)\,dt+\sqrt{\gamma U_t} \,dW_t,
\]
where the invariant distribution is known to be the Gamma distribution
$\Gamma_{2c/\gamma,2cb/\gamma}$, here, the explicit form of the
density is
unknown. It is known (see, e.g., \cite{IkedaWatanabe1989}, Example
IV.8.2, page 237) that $U$ hits $0$ if and only if
$2cb/\gamma<1$. Hence, we may expect that $Y=Y^{\gamma,c,\theta}$
hits $E$ only
at $((2c\theta_2/\gamma,\infty)\times\{0\})\cup(\{0\}\times
(2c\theta_1/\gamma,\infty))$.
Compare this with the fact that $0\in E$ is not hit by the infinite
$\gamma$ process $X^{c,\theta}$ (see Theorem \ref{T1.6}).

\section[Convergence as $\gamma\to\infty$: Proofs of
Theorems 1.3, 1.4]{Convergence as $\gamma\to\infty$: Proofs of
Theorems \protect\ref{T1.3}, \protect\ref{T1.4}}
\label{A3}

\subsection{Construction of the process}
\label{A3.1}

Recall the definitions of $p_t$, $\mathcal{S}$ and $X^{c,\theta}$ in
Definition \ref{D1.1}. In order for the definition to make sense, we
still have to show, in Proposition \ref{P3.2} below, that $p_t$ is
indeed a Markov kernel and that the Chapman--Kolmogorov equation holds.
We prepare for Proposition \ref{P3.2} with a lemma.

Recall the definitions of $C$, $\Xi$, $D$ and ${\mathcal{F}}^D$ in
(\ref{E1.8}),
(\ref{E1.8b}) and (\ref{E1.9}).
\begin{lemma}
\label{L3.1}
\textup{(i)}
$D$ has the Markov property, that is, for $x,y\in[0,\infty)^2$ and
$A\subset E$ measurable, we have
\[
\mathbf{P}[D_{x+y}\in A\mid{\mathcal{F}}_x^D]=Q_{y+D_x}(A).
\]

\smallskipamount=0pt
\begin{longlist}[(ii)]
\item[(ii)]
For $f\dvtx E\to{\mathbb{C}}$ bounded and measurable and $r\geq0$, we have
\[
\int_E f(rz) Q_x(dz)=\int_E f \,dQ_{rx}.
\]
\item[(iii)]
Furthermore,
\[
\int_E Q_x(dz) Q_{rz+y} = Q_{rx+y}.
\]
\end{longlist}
\end{lemma}
\begin{pf}
(i)
Let ${\mathcal{F}}^B$ denote the filtration generated by the Brownian
motion $B$
and let ${\mathcal{F}}^B_{\tau_x}$ denote the $\sigma$-algebra of
the $\tau
_x$ past of $B$ [recall (\ref{E1.7})]. Note that ${\mathcal
{F}}^B_{\tau
_x}\supset{\mathcal{F}}^D_x$.

For $x'\in[0,\infty)^2$, denote by $\mathbf{P}_{-x'}$ the law of $B$ when
started at $B_0=-x'$. Hence, by spatial homogeneity, for $x'\leq x$, we have
\[
\mathbf{P}_{-x'} [B_{\tau_{x+y}}+(x+y)\in A ]=Q_{y+(x-x')}(A).
\]
Choosing $x'=-B_{\tau_x}$, we infer that
\[
\mathbf{P}_{B_{\tau_x}} [B_{\tau_{x+y}}+(x+y)\in A ]=Q_{y+D_x}(A).
\]
We now apply the strong Markov property of $B$ to obtain
\begin{eqnarray*}
\mathbf{P}[D_{x+y}\in A\mid{\mathcal{F}}^D_x]
&=&\mathbf{E} \bigl[\mathbf{P}_0[B_{\tau_{x+y}}+(x+y)\in A\mid
{\mathcal{F}}^B_{\tau_x}]\mid{\mathcal{F}}^D_x \bigr]\\
&=&\mathbf{E} \bigl[\mathbf{P}_{B_{\tau_x}}[B_{\tau_{x+y}}+(x+y)\in
A]\mid{\mathcal{F}}
^D_x \bigr]\\
&=&\mathbf{E} [Q_{y+D_x}(A)\mid{\mathcal{F}}^D_x ]=Q_{y+D_x}(A).
\end{eqnarray*}

\smallskipamount=0pt
\begin{longlist}[(iii)]
\item[(ii)]
This follows from the spatial homogeneity of $B$.

\item[(iii)]
Recall that $D_{rx}$ has distribution $Q_{rx}$.
Hence, by (ii) and (i), we get\end{longlist}
\begin{eqnarray*}
\int_E Q_x(dz) Q_{rz+y}(A)&=&\int_E Q_{rx}(dz) Q_{z+y}(A)\\
&=&\mathbf{E}[Q_{y+D_{rx}}(A)]=\mathbf{P}[D_{rx+y}\in A]\\
&=&Q_{rx+y}(A).
\end{eqnarray*}
\upqed\end{pf}
\begin{proposition}
\label{P3.2}
$(\mathcal{S}_t)_{t\geq0}$ defined in Definition \ref{D1.1} is a
Markov semigroup.
\end{proposition}
\begin{pf}
Recall that $x\mapsto Q_x$ is a continuous map. Hence, for open sets
$A\subset E$, the map $x\mapsto Q_x(A)$ is lower semicontinuous, by the
portmanteau theorem (see, e.g., \cite{Klenke2008e}, Theorem 13.16),
and is hence measurable. Hence, $x\mapsto Q_x(A)$ is measurable for all
Borel sets $A\subset E$. It remains to check the Chapman--Kolmogorov
equation for $(p_t)$.
By Lemma \ref{L3.1}(iii), we infer
that
\begin{eqnarray*}
\int_Ep_t(x,dy) p_s(y,\bolds{\cdot})
&=&\int_EQ_{e^{-ct}x+(1-e^{-ct})\theta}(dy)
Q_{e^{-cs}y+(1-e^{-cs})\theta}\\
&=&Q_{e^{-c(t+s)}x+e^{-cs}(1-e^{-ct})\theta+(1-e^{-cs})\theta}\\
&=&Q_{e^{-c(t+s)}x+(1-e^{-c(t+s)})\theta}\\
&=&p_{t+s}(x,\bolds{\cdot}).
\end{eqnarray*}
\upqed\end{pf}

\subsection[Duality and proof of finite-dimensional distributions
convergence (Theorem 1.3)]{Duality and proof of finite-dimensional distributions
convergence (Theorem~\protect\ref{T1.3})}
\label{A3.2}

In this section, we prove the convergence of the finite-dimensional
distributions of $Y^{\gamma,c,\theta}$ to those of $X=X^{c,\theta}$
by means of a duality relation. For $Y^{\gamma,c,\theta}$, we have
already established the duality,\vspace*{1pt} in Proposition \ref{P2.1}.
We now come to the duality for $X$. Recall the definition of $\tilde Y$ from
(\ref{E2.11}). We will need as initial values only
$\tilde y\in E\times[0,\infty)^2$. Note that, in this case, the
process $Z$ is constant in time and the process $\tilde Y$ is given by
the deterministic equation
%
%
\begin{equation}
\label{E3.1}
\tilde Y_t= \bigl(e^{-ct}\tilde y(1), (1-e^{-ct})\tilde y(1)+\tilde
y(2) \bigr).
\end{equation}
Hence, $\tilde Y$ can be understood as a deterministic Markov process
with state space $E\times[0,\infty)^2$. Recall $H$ from (\ref{E2.12})
and $F$ from (\ref{E2.2}).
\begin{proposition}
\label{P3.4}
$X$ and $\tilde Y$ are dual in the sense that for all initial
conditions $X_0=x\in E$, $\tilde Y_0=\tilde y\in E\times[0,\infty)^2$
and for all $t\geq0$, we have
%
%
\begin{equation}
\label{E3.2}
\mathbf{E}_x [H ((X_t,\theta),\tilde y ) ]
=\mathbf{E}_{\tilde y} [H ((x,\theta),\tilde Y_t ) ].
\end{equation}
In particular, we get
%
%
\begin{equation}
\label{E3.3}
\mathbf{E}_x[F(X_t,z)]=F (x,e^{-ct}z ) F \bigl(\theta,
(1-e^{-ct} )z \bigr)\qquad
 \mbox{for }x\in[0,\infty)^2,z\in E,\hspace*{-32pt}
\end{equation}
and the distribution of $X_t$ is determined by (\ref{E3.3}).
\end{proposition}
\begin{pf}
As $\tilde Y$ is deterministic, (\ref{E3.2}) and (\ref{E3.3}) are
equivalent and so we only need to show (\ref{E3.3}). Since $z\in E$, by
Corollary \ref{C2.3}(iii), the left-hand side of (\ref{E3.3}) equals
\begin{eqnarray*}
\int_E F(y,z) Q_{e^{-ct}x+(1-e^{-ct})\theta}(dy)
&=&F \bigl(e^{-ct}x+(1-e^{-ct})\theta,z \bigr)\\
&=&F (x,e^{-ct}z ) F \bigl(\theta,(1-e^{-ct})z \bigr).
\end{eqnarray*}
By Corollary \ref{C2.4}, equation (\ref{E3.3}) determines the
distribution of $X_t$.
\end{pf}

We are now ready to prove Theorem \ref{T1.3}.
\begin{pf*}{Proof of Theorem \ref{T1.3}}
As both $X^{c,\theta}$ and $Y^{\gamma,c,\theta}$ are Markov processes,
it is easy to see that for convergence of finite-dimensional
distributions, it is enough to show that for any $t\geq0$, $x\in E$
and $(x_{\gamma})_{\gamma\geq0}$ in $[0,\infty)^2$ such that $\lim
_{\gamma\rightarrow\infty} x_{\gamma}\rightarrow x$, we have
%
%
\begin{equation}
\label{E3.4}
{\mathcal{L}}_{x_{\gamma}} [Y^{\gamma,c,\theta}_t
]\stackrel{\gamma\to\infty}{\longrightarrow}{\mathcal{L}}_x
[X^{c,\theta}_t ] \qquad\mbox{weakly. }
\end{equation}
As shown in the proof of Corollary \ref{C2.4}(i), $\mathcal{D}_0$ is dense
in $C_l([0,\infty)^2)$. Hence, it is enough to consider $F(\bolds
{\cdot},z)$,
$z\in[0,\infty)^2$, as test functions. Denote by
$Z^\gamma$ the process defined in (\ref{E2.10}) started at $Z^\gamma
_0=z$. For
$\gamma=1$, we drop the superscript, that is,
$Z:=Z^1$. Denote by $Z_\infty$ the almost sure limit of $Z_t$ as $t\to
\infty$
and recall that its distribution is $Q_z$. Note that, due to Brownian scaling,
$(Z^\gamma_t)_{t\geq0}\stackrel{\mathcal{D}}{=}(Z_{\gamma
t})_{t\geq0}$. Hence, by
Proposition \ref{P2.1}, we have
\begin{eqnarray*}
\mathbf{E}_{x_{\gamma}} [F (Y^{\gamma,c,\theta}_t,z
) ]
&=&
\mathbf{E}_z \biggl[F (x_{\gamma}, e^{-ct}Z_{\gamma t} )
F \biggl(\theta, \int_0^t c  e^{-cr}Z_{\gamma r} \,dr \biggr) \biggr]\\
&\stackrel{\gamma\to\infty}{\longrightarrow}&
\mathbf{E}_z \bigl[F (x, e^{-ct}Z_{\infty} )
F \bigl(\theta, (1 - e^{-ct})Z_{\infty} \bigr) \bigr]\\
&=&\int_E F (x,e^{-ct}y )F \bigl(\theta, (1-e^{-ct})y \bigr)
Q_z(dy)\\
&=&\int_E \mathbf{E}_x[F(X_t,y)] Q_z(dy)\\
&=&\mathbf{E}_x \biggl[\int_EF(X_t,y) Q_z(dy) \biggr]\\
&=&\mathbf{E}_x [F(X_t,z) ],
\end{eqnarray*}
where the fourth line follows by (\ref{E3.3}) and
the last equality follows by Corollary~\ref{C2.3}(iii).
\end{pf*}
\begin{remark}
\label{R3.5}
We could also define $X^{c,\theta}$ in Definition \ref{D1.1} for initial
values
$x\in[0,\infty)^2$ (instead of $E$ only). This means that
$X^{c,\theta}$ starts life with a jump from
$x$ to a random point on $E$ chosen according to $Q_x$ and then
continues with the usual dynamics. Clearly, this process does not have
a c{\`a}dl{\`a}g version (due
to the jump at time $0$) and its transition semigroup is not strongly
continuous at $0$. Nevertheless, the proof of Theorem \ref{T1.3} shows that
that theorem also holds for this process and hence for
$Y^{\gamma,c,\theta}_0=X^{c,\theta}_0=x\in[0,\infty)^2$.
\end{remark}

\subsection[Proof of the $L^p$-convergence (Theorem 1.4)]{Proof of the
$L^p$-convergence (Theorem \protect\ref{T1.4})}
\label{A3.3}

We prepare for the proof of Theorem \ref{T1.4} with two lemmas.
\begin{lemma}
\label{L3.6}
Let $B=(B_1,B_2)$ be a planar Brownian motion started at
$(B_{1,0},B_{2,0})=(u,v)\in[0,\infty)^2$ and let
\[
\tau=\inf \{t>0\dvtx B_t\notin(0,\infty)^2 \}.
\]
Then, for any $p\in[1,2)$, we have
\[
\mathbf{E} [\tau^{p/2} ]\leq
\frac{2}{2-p} \biggl(\frac{2}{\pi} \biggr)^{p/2}(uv)^{p/2} < \infty.
\]
\end{lemma}

More generally, one could show for the exit time of a cone with angle
$2\alpha$ (here, $\alpha=\pi/4$) that $\mathbf{E}[\tau
^{p/2}]<\infty$ if
and only if $p\alpha<\pi/2$ (see \cite{Burkholder1977}, equation
(3.8)). We give the short proof here in order to be
self-contained.
\begin{pf}
By the reflection principle and independence of $B_1$ and $B_2$, we get
\[
\mathbf{P}[\tau>t]=4\mathcal{N}_{0,t}(0,u)\mathcal{N}_{0,t}(0,v),
\]
where $\mathcal{N}_{0,t}(a,b)=(2\pi t)^{-1/2}\int_a^be^{-r^2/2t}\,dr$
is the
centred normal distribution with variance $t$. Hence,
\begin{eqnarray*}
\mathbf{E} [\tau^{p/2} ]&=&\int_0^\infty\mathbf{P}
[\tau>t^{2/p}
] \,dt\\
&\leq&\int_0^\infty1\wedge \biggl(\frac2\pi  uv t^{-2/p}  \biggr) \,dt
=\frac{2}{2-p} \biggl(\frac{2}{\pi} \biggr)^{p/2}(uv)^{p/2}.
\end{eqnarray*}
\upqed\end{pf}
\begin{lemma}
\label{L3.7}
For every $(u,v)\in[0,\infty)^2$, every $p\in[1,2)$ and every
$i=1,2$, we have
\[
\int_E x_i^{p} Q_{(u,v)}(dx)
\leq |u^2-v^2|^{p/2} + \frac{2^{p/2}(uv)^{p/2}}{\cos(p\pi
/4)}<\infty.
\]
\end{lemma}
\begin{pf}
This can be verified by means of an explicit computation using the
density formula of $Q_{(u,v)}$ in (\ref{E2.5}).
\end{pf}

Note that finiteness of the expression on the left-hand side in
Lemma \ref{L3.7} (which is what we need in the proof of Theorem \ref
{T1.4}) could also be inferred without computations by the
Burkholder--Davis--Gundy inequality and Lemma \ref{L3.6}.
\begin{pf*}{Proof of Theorem \ref{T1.4}}
(i) By Lemma \ref{L3.7}, we have
\[
\mathbf{E} [(X^{c,\theta}_{i,t})^p ]=\int_E y_i^p
Q_{e^{-ct}x+(1-e^{-ct})\theta}(dy) <\infty.
\]

Fix $t>0$ and define
\[
M^{t}_{i,s}:=e^{-ct}x_i+(1-e^{-ct})\theta_i +\int_0^se^{c(r-t)}\sqrt
{\gamma Y^{\gamma,c,\theta}_{1,r}Y^{\gamma,c,\theta}_{2,r}} \,dW_{i,r}.
\]
Let $\langle M^t_1\rangle=\langle M^t_2\rangle$ denote the square
variation process of both $M^t_1$ and $M^t_2$.
Note that
$M^{t}_{i,t}
=Y^{\gamma,c,\theta}_{i,t}\geq0$ and that $M^t_i$ is a martingale
and thus
%
%
\begin{equation}
\label{E3.5}
M^t_{i,s}=\mathbf{E}[M^t_{i,t}\mid M^t_{i,s}]\geq
0 \qquad\mbox{for all }s\in[0,t].
\end{equation}

Now, $(M^{t}_s)_{s\geq0}$ is an isotropic diffusion in ${\mathbb
{R}}^2$ and is
hence a time-transformed planar Brownian motion. That is, there exists
a planar Brownian motion $B$ (with respect to some right-continuous
complete filtration ${\mathcal{F}}$) started at $B_0=
e^{-ct}x+(1-e^{-ct})\theta
$ such that each $\langle M^t_{1,\bolds{\cdot}}\rangle_s$ is an
${\mathcal{F}}$
stopping time and such that $B_{\langle M^t_{1,\bolds{\cdot}}\rangle
_s}=M^t_{s}$ for all $s\geq0$.

Define the ${\mathcal{F}}$ stopping times
\[
\tau:=\inf \{s>0\dvtx B_s\notin(0,\infty) ^2 \}
\quad\mbox{and}\quad
\tau_0:=\inf \{s>0\dvtx B_s\notin[0,\infty) ^2 \}.
\]
Clearly, we have $\tau=\tau_0$ almost surely and, hence, by (\ref{E3.5}),
\[
\langle M^t_{1,\bolds{\cdot}}\rangle_t\leq\tau_0=\tau \qquad\mbox{a.s.}
\]
Using the Burkholder--Davis--Gundy inequality for the martingale
$(B_{i,s})_{s\geq0}$ yields (see Lemma \ref{L3.6})
\[
\mathbf{E} \Bigl[\sup_{s\leq\tau}B_{i,s}^p \Bigr]\leq2^{p-1}
\bigl(B_{i,0}^p+(4p)^p \mathbf{E}[\tau^{p/2}] \bigr)<\infty.
\]
Hence, $(|B_{i,\tau\wedge s}|^p)_{s\geq0}$ is uniformly integrable
and we can apply the optional sampling theorem to obtain
\[
\mathbf{E} [(Y^{\gamma,c,\theta}_{i,t})^p ]=\mathbf{E}
\bigl[\bigl(B_{i,\langle
M^t_{1,\bolds{\cdot}}\rangle_t}\bigr)^p \bigr]
\leq\mathbf{E} [(B_{i,\tau})^p ]=\mathbf{E}
[(X^{c,\theta
}_{i,t})^p ].
\]

(ii) By Theorem \ref{T1.3} and the Skorohod embedding
theorem, we may construct all processes on one probability space such
that $Y^{\gamma,c,\theta}_t\to X^{c,\theta}_t$ almost surely as
$\gamma\to\infty$. By part (i), the $p$th moments of $Y^{\gamma
,c,\theta}_{i,t}$, $\gamma\geq0$, are uniformly integrable and so we
have the desired $L^p$-convergence.
\end{pf*}

\section[The strong construction (proofs of Theorems 1.6 and 1.7)]{The strong construction (proofs
of Theorems \protect\ref{T1.5} and \protect\ref{T1.6})}
\label{A4}

Recall the definitions of $C$, $\Xi$ and $D$ in (\ref{E1.8}) and
(\ref{E1.9}).
\begin{lemma}
\label{L4.1}
The map $x\mapsto D_x$ is c{\`a}dl{\`a}g.
\end{lemma}
\begin{pf}
This follows from continuity of $B$ and the definition of $\tau_x$.
\end{pf}
\begin{pf*}{Proof of Theorem \ref{T1.5}}
From Lemmas \ref{L3.1} and \ref{L4.1}, we infer that $X^{\bar c,\bar
\theta}$ has the Markov property and c{\`a}dl{\`a}g paths. It remains
to show (\ref{E1.10}).

By Lemma \ref{L3.1}, for $x,z\in E$, $A\subset E$ measurable and
$0\leq s<t$, we have (with $\mathbf{P}_x$ denoting the probability law of
$X^{\bar c,\bar\theta}_t$, as defined in Theorem \ref{T1.5})
\begin{eqnarray*}
p_{s,t}(z,A)
&=&\mathbf{P}_x [X^{\bar c,\bar\theta}_t\in A\mid X^{\bar c,\bar
\theta
}_s=z ]\\
&=&\mathbf{P} \bigl[C(0,t) D_{x+\Xi(0,t)}\in A\mid
D_{x+\Xi(0,s)}=C(0,s)^{-1}z \bigr]\\
&=&Q_{C(0,s)^{-1}z+\Xi(s,t)} (C(0,t)^{-1}A )\\
&=&Q_{C(s,t)z+C(0,t)\Xi(s,t)} (A ).
\end{eqnarray*}
\upqed\end{pf*}
\begin{pf*}{Proof of Theorem \ref{T1.6}}
If $c\theta=0$, then $X^{c,\theta}$ is the deterministic process
$X^{c,\theta}_t=e^{-ct}X^{c,\theta}_0$ and hence $0$ is polar.

Now, assume that $c\theta\neq0$.
Le Gall and Meyre \cite{LeGallMeyre1992} show that almost surely, for
all $z\in(0,\infty)^2$,
the planar Brownian motion $B$ does not leave the cone $[-z,\infty)$
first at $-z$. More formally, consider the event
\[
A:= \{B_{\tau_z}\neq-z\mbox{ for all }z\in(0,\infty)^2 \}.
\]
Theorem 1 of \cite{LeGallMeyre1992} then implies that $\mathbf
{P}[A]=1$ (in
fact, they show that \textit{no} rectangular cone is first left at its
vertex, not only northeast cones $[z,\infty)$). Now, by (\ref{E1.11}),
we have
\begin{eqnarray*}
\{X^{c,\theta}_t\neq0\mbox{ for all }t>0 \}
&=& \{D_{x+r\theta}\neq0\mbox{ for all }r>0 \}\\
&=& \{B_{\tau_{x+r\theta}}\neq x+r\theta\mbox{ for all }r>0
\}
\supset A.
\end{eqnarray*}
This shows the claim of Theorem \ref{T1.6}.
\end{pf*}
%
\section{The infinitesimal dynamics of $X^{c,\theta}$}
\label{A5}

In this section, we give a description and construction of the infinite
rate mutually catalytic branching process $X$ in terms of its
infinitesimal characteristics. To this end, we will define a linear
operator ${\mathcal{G}}^{c,\theta}$ that:
\begin{longlist}
\item
defines the contraction semigroup of $X$ in the sense of the
Hille--Yoshida theorem (Theorem \ref{T5.3});
\item
defines a well-posed martingale problem whose unique solution is $X$
(Theorem \ref{T5.4}).
\end{longlist}
%
\subsection{Results}
\label{A5.1}

Recall, from Definition \ref{D1.1}, that the linear operator $\mathcal{S}_t$
on $C_l(E)$ is defined by
\[
\mathcal{S}_tf(x):=\int_E f(y) p_t(x,dy)=\int_E f(y)
Q_{e^{-ct}x+(1-e^{-ct})\theta}(x,dy).
\]
In order to define the generator of $\mathcal{S}=(\mathcal
{S}_t)_{t\geq0}$, we will
need to study (for suitable functions $f$) the limit
%
%
\begin{equation}
\label{E5.1}
\lim_{t\downarrow0}t^{-1}\bigl(\mathcal{S}_tf(x)-f(x)\bigr)=\lim
_{{\varepsilon}\to0}{\varepsilon}
^{-1} \biggl(\int f \,dQ_{x+\varepsilon c(\theta-x)}-f(x) \biggr).
\end{equation}
In the sequel, we will use the shorthand notation
\[
\partial_1 f(u,v):=\frac{\partial}{\partial u}f(u,v)
\quad\mbox{and}\quad
\partial_2 f(u,v):=\frac{\partial}{\partial v}f(u,v).
\]
In order to define what we mean by a suitable function, we introduce
the subspace $C_l^2(E)\subset C_l(E)$.
\begin{definition}
\label{D5.1}
Let $C_l^2(E)\subset C_l(E)$ be the subspace of such functions $f\in C_l(E)$:
\begin{longlist}
\item
whose partial derivatives
$\partial_1f$ and $\partial_2f$ exist on $(0,\infty)\times\{0\}$ and
$\{0\}\times(0,\infty)$, respectively, are continuous, can be
continuously extended to $\{0\}\times[0,\infty)$ and fulfill
%
%
\begin{equation}
\label{E5.2}
\lim_{u\to\infty}u\,\partial_1f(u,0)=\lim_{v\to\infty}v\,\partial_2f(0,v)=0;
\end{equation}
\item
whose partial second derivatives $\partial_1^2f$ and $\partial_2^2f$
exist on $(0,\infty)\times\{0\}$ and $\{0\}\times(0,\infty)$,
respectively, and are such that
%
%
\begin{equation}
\label{E5.3}
\|f\|_{2,\infty}:=
\sup_{r\in(0,\infty)}r \bigl(|\partial_1^2f(r,0)|+|\partial
^2_2f(0,r)| \bigr)<\infty.
\end{equation}
\end{longlist}
\end{definition}

Note that, for $f\in C_l^2(E)$, we have
%
%
\begin{equation}
\label{E5.4}
\|f\|_{1,\infty}:=\sup_{r\in[0,\infty)}
\bigl(|\partial_1f(r,0)|+|\partial_2f(0,r)| \bigr)<\infty.
\end{equation}

In order to get an explicit formula for the limit in (\ref{E5.1}), we
define the vague limits (for $u,v>0$)
\[
\nu_{(0,v)}:=\mathop{v\mbox{-}\mathrm
{lim}}_{{\varepsilon}\downarrow0} {\varepsilon
}^{-1}Q_{({\varepsilon},v)}
\quad\mbox{and}\quad
\nu_{(u,0)}:=\mathop{v\mbox{-}\mathrm
{lim}}_{{\varepsilon}\downarrow0}{\varepsilon
}^{-1}Q_{(u,{\varepsilon})}.
\]
$\nu_{(u,0)}$ can be thought of as the ``L{\'e}vy measure'' of the
next jump when the actual position is $(u,0)$ and the drift is in
direction of $(0,1)$.
In order to formalize this, for the drift in direction $(0,1)$, we
define the linear operator ${\mathcal{G}}_2$ on $C_l^2(E)$ by
${\mathcal{G}}_2
f(x)=\partial_2 f(x)$ if
$x_1=0$ and
\[
{\mathcal{G}}_2 f(x)=\int_E [f(y)-f(x)-(y_1-x_1)\,\partial_1f(x)] \nu_{x}(dy)
\qquad \mbox{if }x_1>0.
\]
For the drift in direction $(1,0)$, we define ${\mathcal{G}}_1$ similarly.
Note that $\nu_x$ is not a finite measure and that the integral of
$y_1-u$ with respect to $\nu_{(u,0)}$ is well defined only as a Cauchy
principal value and, as such, equals zero. Hence, this term in the
integral is needed in order for the integral to be well defined in the
usual sense.
In Lem\-ma~\ref{L5.5} below, we will show that ${\mathcal{G}}_1f$ and
${\mathcal{G}}_2f$
are, in fact, well defined and are in $C_l(E)$.

Due to spatial homogeneity of planar Brownian motion, we have a scaling
relation that helps to get rid of the many different $\nu_x$ in the
definition of ${\mathcal{G}}_1$ and~${\mathcal{G}}_2$:
\[
\int_E f(x) \nu_{(u,0)}(dx) = \frac{1}{u}\int_E f(ux) \nu_{(1,0)}(dx).
\]
Furthermore, letting $f^\dagger((x_1,x_2)):=f((x_2,x_1))$, by
symmetry, we have
\[
\int_E f(x) \nu_{(0,v)}(dx) = \int_E f^\dagger(x) \nu_{(v,0)}(dx)
  = \frac{1}{v}\int_E f^\dagger(vx) \nu_{(1,0)}(dx).
\]

Hence, we can express ${\mathcal{G}}_1$ and ${\mathcal{G}}_2$ in
terms of
%
%
\begin{equation}
\label{E5.5}
\nu:=\nu_{(1,0)}.
\end{equation}
Using the explicit form of the density of $Q_{(1,{\varepsilon})}$ in
(\ref{E2.5})
and letting ${\varepsilon}\to0$, we get that the $\sigma$-finite
measure $\nu$
on $E$ has a one-dimensional Lebesgue density given by
%
%
\begin{equation}
\label{E5.6}
\nu (d(u,v) )=\cases{
\displaystyle\frac4\pi \frac{u}{(1-u)^2 (1+u)^2}\, du,&\quad
if $v=0$,\cr
\displaystyle\frac4\pi \frac{v}{ (1+v^2 )^2}
\,dv,&\quad if $u=0$.}
\end{equation}
${\mathcal{G}}_1$ and ${\mathcal{G}}_2$ can now be written as
%
%
\begin{equation}
\label{E5.7}
{\mathcal{G}}_2f (x)=\cases{
\displaystyle\partial_ 2f(x),&\quad if $x_1=0$,\vspace*{1pt}\cr
\displaystyle\frac1{x_1}\int_E [f(x_1y)-f(x)\vspace*{1pt}\cr
\hspace*{28.4pt}{}
\displaystyle-x_1(y_1-1) \,\partial_1 f(x) ]
\nu(dy),&\quad if $x_1>0$,}
\end{equation}
and
%
%
\begin{equation}
\label{E5.8}
{\mathcal{G}}_1f=({\mathcal{G}}_2f^\dagger)^\dagger.
\end{equation}
Finally, we define the operator ${\mathcal{G}}^{c,\theta}$ on
$C_l^2(E)$ with
domain $\mathcal{D}({\mathcal{G}}^{c,\theta})=C_l^2(E)$ that
determines the
infinitesimal characteristics of the process $X=X^{c,\theta}$:
%
%
\begin{equation}
\label{E5.9}
{\mathcal{G}}^{c,\theta}f(x)=\sum_{i=1}^2 c(\theta_i-x_i){\mathcal
{G}}_if(x).
\end{equation}
\begin{lemma}
\label{L5.2}
The operator ${\mathcal{G}}^{c,\theta}$ is well defined. That is, for
$f\in
C_l^2(E)$, the expressions in (\ref{E5.9}) and (\ref{E5.7}) are well
defined and we have ${\mathcal{G}}^{c,\theta}f\in C_l(E)$.
\end{lemma}

This lemma will be proven in Section \ref{A5.2}.
\begin{theorem}[($X^{c,\theta}$ via its generator)]
\label{T5.3}
\textup{(i)}
For every $f\in C_l^2(E)$, we have, pointwise, for all $x\in E$,
that
%
%
\begin{equation}
\label{E5.10}
{\mathcal{G}}^{c,\theta}f(x)=\lim_{{\varepsilon}\downarrow
0}{\varepsilon}^{-1} \biggl(\int_E
f \,dQ_{x+{\varepsilon} c(\theta-x)} -f(x) \biggr)
=\lim_{t\downarrow0}\frac{\mathcal{S}_t f(x)-f(x)}{t}.\hspace*{-35pt}
\end{equation}

\smallskipamount=0pt
\begin{longlist}[(ii)]
\item[(ii)]
The operator ${\mathcal{G}}^{c,\theta}$ on $C_l(E)$ is closable and
its closure
generates the contraction semigroup $\mathcal{S}$ of the process
$X^{c,\theta}$.
\end{longlist}
\end{theorem}

The theorem will be proven in Section \ref{A5.2} using the classical
Hille--Yoshida theorem.

A different, and more modern, approach to constructing Markov processes
from their infinitesimal dynamics is the martingale problem technique
due to Stroock and Varadhan.

Recall from (\ref{E2.17}) that $V\subset C_l^2(E)$ is the vector space
spanned by $\{F(\bolds{\cdot},z)$, $z\in E\}$. Define the linear
operator ${\mathcal{G}}
^{c,\theta}$ on $V$ by (\ref{E5.9}) and (\ref{E5.7}). By
Theorem \ref{T5.3}(i), we obtain for $z\in E$ [using Corollary \ref{C2.3}(iii) in
the second line and Lemma \ref{L2.2} in the last line] that
%
%
\begin{eqnarray}
\label{E5.11}\qquad
{\mathcal{G}}^{c,\theta} F(\bolds{\cdot},z)(x)
&=&\lim_{{\varepsilon}\downarrow0}{\varepsilon}^{-1} \biggl(\int_E
F(y,z) \,dQ_{x+{\varepsilon} c(\theta-x)}(dy) -F(x,z) \biggr)\nonumber\\
&=&\lim_{{\varepsilon}\downarrow0}{\varepsilon}^{-1} \bigl(
F \bigl(x+{\varepsilon} c(\theta-x),z \bigr)-F (x,z ) \bigr)\nonumber\\[-8pt]\\[-8pt]
&=&\langle c (\theta-x), \nabla_1 F(x,z) \rangle
\nonumber\\
&=&F (x,y ) [c(\theta-x)\diamond z].\nonumber
\end{eqnarray}

Hence, (\ref{E5.11}) is enough to define ${\mathcal{G}}^{c,\theta}$
on $V$ and
we do not really need the measure $\nu$ from (\ref{E5.7}) here.

A solution of the $({\mathcal{G}}^{c,\theta},V)$ martingale problem
is an
$E$-valued measurable stochastic process $X$ such that
\[
M_t:=F(X_t,z)-\int_0^t \bigl(c(\theta-X_s)\diamond z\bigr) F(X_s,z) \,ds
\]
is a (${\mathbb{C}}$-valued) martingale. A martingale problem is said
to be
\textit{well posed} if, for every probability measure $\mu$ on $E$,
there exists a solution $X$ with ${\mathcal{L}}[X_0]=\mu$ (existence)
and any
two solutions have the same finite-dimensional distributions
(uniqueness). In this case, $X$ is a Markov process (see \cite{EthierKurtz1986},
Theorem 4.4.2(a)).
\begin{theorem}[(Martingale problem characterization of $X^{c,\theta}$)]
\label{T5.4}
The martingale problem $({\mathcal{G}}^{c,\theta},V)$ is well posed
and its
unique solution is $X^{c,\theta}$.
\end{theorem}

This theorem will be proven in Section \ref{A5.3}.

\subsection[The Hille--Yoshida approach (proof of Theorem 5.3)]{The Hille--Yoshida
approach (proof of Theorem \protect\ref{T5.3})}
\label{A5.2}

Lemma \ref{L5.2} and part~(i) of Theorem \ref{T5.3} are direct
consequences of the following two lemmas.

Let $e_1=(1,0)$ and $e_2=(0,1)$.
\begin{lemma}
\label{L5.5}
For $f\in C_l^2(E)$, $x\in E$ and $i=1,2$, the expression ${\mathcal{G}}_if(x)$
from (\ref{E5.7}) and (\ref{E5.8}) is well defined and we have
%
%
\begin{equation}
\label{E5.12}
\lim_{{\varepsilon}\downarrow0} {\varepsilon}^{-1} \biggl(\int_E
f \,dQ_{x+{\varepsilon}
e_i} -f(x) \biggr)={\mathcal{G}}_if(x).
\end{equation}
\end{lemma}
\begin{lemma}
\label{L5.6}
For $f\in C_l^2(E)$, we have ${\mathcal{G}}^{c,\theta}f\in C_l(E)$.
\end{lemma}
\begin{pf*}{Proof of Lemma \ref{L5.5}}
For $x=(0,0)$, since $Q_{{\varepsilon}e_i}=\delta_{{\varepsilon
}e_i}$, this is the very
definition of ${\mathcal{G}}_i$. For $u\neq(0,0)$,
by linear scaling and symmetry, it is enough to consider the case
$x=(1,0)$. If $i=1$, then the left-hand side of (\ref{E5.12})
equals
\[
\lim_{{\varepsilon}\downarrow0} {\varepsilon}^{-1}
\bigl(f(1+{\varepsilon},0)-f(1,0)
\bigr)=\partial_1
f(1,0)=({\mathcal{G}}_1f)(1,0).
\]
Now, consider $i=2$. It is a simple exercise to compute that for
every ${\varepsilon}>0$,
\begin{eqnarray*}
\frac4\pi\int_0^\infty\frac{r(r-1)}{4{\varepsilon
}^2+(r^2+{\varepsilon}^2-1)^2} \,dr
&=&\frac{2}{\pi} {\varepsilon}^{-1}\arctan({\varepsilon})\\
&=&\frac4\pi\int_0^\infty\frac{s}{4{\varepsilon
}^2+(s^2-{\varepsilon}^2+1)^2} \,ds.
\end{eqnarray*}
Hence, if we let $g(y):=(y_1-1)\,\partial_1f(1,0)$, then, for every
${\varepsilon}>0$,
\[
\int_E \bigl(g(y)-g(1,0) \bigr) Q_{(1,{\varepsilon})} (dy )
= 0.
\]
Hence, we can replace $f$ by $f-g$. Now, $f-g$ is twice differentiable,
has at most linear growth and
$\partial_1 (f-g)(1,0)=0$. Hence,
\[
\sup_{u\geq0,  u\neq1}\frac{|(f-g)(u,0)-f(1,0)|}{(u-1)^2}<\infty.
\]
This allows us to use dominated
convergence in the following computation to obtain
\begin{eqnarray*}
&&\lim_{{\varepsilon}\downarrow0} {\varepsilon}^{-1} \biggl(\int
_E f \,dQ_{(1,0)+{\varepsilon}
e_2} -f(1,0) \biggr)\\
&&\qquad=\lim_{{\varepsilon}\downarrow0}{\varepsilon}^{-1}\int
[f(x)-f(1,0)-(x_1-1)\,\partial_1f(1,0) ]Q_{(1,{\varepsilon})}(dx)\\
&&\qquad=\lim_{{\varepsilon}\downarrow0} \biggl(
\frac4\pi\int_0^\infty
\frac{u
[(f-g)(u,0)-f(1,0)]}{4{\varepsilon}^2+ (u^2+{\varepsilon
}^2-1 )^2} \,du\\
&&\qquad\quad\hspace*{19.9pt}{} + \frac4\pi\int_0^\infty
\frac{v [(f-g)(0,v)-f(1,0)]}{4{\varepsilon}^2+
(v^2-{\varepsilon}^2+1 )^2} \,dv
 \biggr)\\
&&\qquad=\lim_{{\varepsilon}\downarrow0} \biggl(
\frac4\pi\int_0^\infty
\frac{u
[(f-g)(u,0)-f(1,0)]}{4{\varepsilon}^2+{\varepsilon
}^4+2{\varepsilon}^2(u+1)(u-1)+(u+1)^2(u-1)^2}
\,du\\
&&\qquad\quad\hspace*{111pt}{}
 + \frac4\pi\int_0^\infty\frac{v
[(f-g)(0,v)-f(1,0)]}{4{\varepsilon}
^2+ (v^2-{\varepsilon}^2+1 )^2} \,dv\biggr)\\
&&\qquad=\frac4\pi\int_0^\infty
\frac{u [(f-g)(u,0)-f(1,0)]}{(u^2-1)^2} \,du\\
&&\qquad\quad{} + \frac4\pi\int_0^\infty
\frac{v [(f-g)(0,v)-f(1,0)]}{(v^2+1)^2} \,dv
\\
&&\qquad=\int_E [f(y)-f(1,0)
-(y_1-1)\,\partial_1 f(1,0) ] \nu(dy)\\
&&\qquad={\mathcal{G}}_2f(1,0).
\end{eqnarray*}
\upqed\end{pf*}
\begin{pf*}{Proof of Lemma \ref{L5.6}}
We have to show that for any $f\in C_l^2(E)$, ${\mathcal{G}}^{c,\theta
}f(x)$ is
continuous in $x\in E$ and has a limit at $\infty$. By (\ref{E5.9}),
it is
enough to derive these properties for $G_i(x):=(\theta_i-x_i){\mathcal{G}}
_if(x)$, $i=1,2$. We
will give the proof only for the case $i=2$ since the case $i=1$ is
analogous.

For $x_1=0$, we have
%
%
\begin{equation}
\label{E5.13}
G_2(x)=G_2(0,x_2)= (\theta_2-x_2)\,\partial_2f(0,x_2).
\end{equation}
This expression is clearly continuous in
$x_2\in[0,\infty)$
and, by (\ref{E5.2}), we have
%
%
\begin{equation}
\label{E5.14}
\lim_{x_2\to\infty}G_2(x)=0.
\end{equation}

Now, consider the case $x_1>0$. Hence, by (\ref{E5.7}),
\[
G_2(x)=\int g(x,y) \nu(dy),
\]
where
\[
g(x,y):=
\frac{\theta_2}{x_1} [f(x_1y)-f(x)
-x_1(y_1-1) \,\partial_1 f(x) ].
\]
Since $f\in C_l^2(E)$, for all $y\in E$, we have:

\begin{longlist}
\item
$x\mapsto g(x,y)$ is continuous on $(0,\infty)\times\{0\}$;
\item
$\lim_{x_1\to\infty}g(x,y)= 0$;
\item
$\lim_{x_1\downarrow0} g(x,y)= \theta_2\,\partial_2f(0,0) y_2$.
\end{longlist}
In order to find an integrable dominating function for $g$, define
$h\dvtx E\to[0,\infty)$ by [recall (\ref{E5.3}) and (\ref{E5.4})]
\[
h(y):=\cases{
\theta_2 \|f\|_{2,\infty} (y_1-1)^2,&\quad
if $y_1\in \bigl(\frac12,\frac{3}{2} \bigr)$,\cr
2\theta_2 \|
f\|_{1,\infty}  (y_1+y_2+1),&\quad otherwise.}
\]
Note that the density of $\nu(dy)$ decays like $1/(y_1+y_2)^3$ as
$y\to\infty$. Furthermore,
\[
(y_1-1)^2\frac{\nu(dy)}{dy_1}=\frac{4}{\pi}\frac{y_1}{(1+y_1)^2}
\]
is bounded on $(1/2,3/2)\times\{0\}$. Hence, we have $\int h \,d\nu
<\infty$.

For all $y\in E$ and
$x_1>0$,
we have
\begin{eqnarray*}
|g(x,y)|&\leq&\frac{\theta_2}{x_1}
\bigl(|f(x_1y)-f(0,0)|+|f(x)-f(0,0)|+x_1(y_1+1)|\partial_1f(x)| \bigr)\\
&\leq&2\theta_2
(y_1+y_2+1) \|f\|_{1,\infty}.
\end{eqnarray*}
Furthermore, recalling (\ref{E5.3}), for $y_1\in(1/2,3/2)$, by Taylor's
formula, we get that
\begin{eqnarray*}
|g(x,y)|&=&\frac{\theta_2}{x_1} \bigl|f\bigl((y_1-1)x+x\bigr)-f(x)
-x_1(y_1-1)\,\partial_1f(x) \bigr|\\
&\leq&\frac{\theta_2}{2}
(y_1-1)^2
\sup_{u\geq x_1/2}x_1 |\partial_1^2f(u,0)|\\
&\leq&\theta_2 \|f\|_{2,\infty} (y_1-1)^2.
\end{eqnarray*}
Hence, in fact, $|g(x,y)|\leq h(y)$ for all
$y\in E$, $x\in(0,\infty)\times\{0\}$
and the dominated convergence theorem yields that $G_2$ shares the
properties (i) and (ii) of $g(x,\bolds{\cdot})$ and that
\[
\lim_{x_1\downarrow0}G_2(x) = \theta_2 \,\partial_2f(0,0)\int y_2
\nu(dy) = \theta_2 \,\partial_2f(0,0) =
G_2(0,0).
\]
Combining this with (\ref{E5.13}) and (\ref{E5.14}), we have $G_2\in C_l(E)$.
\end{pf*}

In order to show part (ii) of Theorem \ref{T5.3}, we will apply the
Hille--Yoshida theorem for generators of contraction semigroups.
Recall, from Corollary \ref{C2.4}, that
$V$ is dense in $C_l(E)$. Also, by Lemma \ref{L2.2}, one can easily
check that
\[
V\subset C_l^2(E).
\]
For each $z\in E$, define the map $u_y\dvtx[0,\infty)\to C_l(E)$ by
$u_y(t):=\mathcal{S}_tF(\bolds{\cdot},y)$.

By \cite{EthierKurtz1986}, Proposition 1.3.4, the operator ${\mathcal{G}}
^{c,\theta}$
on $C_l(E)$ with domain $\mathcal{D}({\mathcal{G}}^{c,\theta
})=C_l^2(E)$ is closable
and its closure
generates (uniquely) the semigroup $(\mathcal{S}_t)_{t\geq0}$ on
$C_l(E)$ if the
following conditions are all fulfilled:
\begin{enumerate}[(a)]
\item[(a)] ${\mathcal{G}}^{c,\theta}$ is dissipative;
\item[(b)]
$u_y(t)\in\mathcal{D}({\mathcal{G}}^{c,\theta})$ for all $t>0$;
\item[(c)]
the map $(0,\infty)\to C_l(E)$, $t\mapsto{\mathcal{G}}^{c,\theta
}u_y(t)$ is
continuous;
\item[(d)] for all $t>0$,
%
%
\begin{equation}
\label{E5.15}
u_y(t)-u_y(0)=\int_0^t{\mathcal{G}}^{c,\theta}u_y(s) \,ds.
\end{equation}
\end{enumerate}
Hence, in order to prove part (ii) of Theorem \ref{T5.3}, it remains
to check (a)--(d).

(a) Let $f\in C_l^2(E)$ and assume that $f$ assumes its
maximum at $x\in E\cup\{\infty\}$. Since $\mathcal{S}_tf(x)\leq
f(x)$ for
all $t\geq0$, equation (\ref{E5.10}) implies that ${\mathcal
{G}}^{c,\theta
}f(x)\leq0$. Hence, ${\mathcal{G}}^{c,\theta}$ fulfills the positive maximum
principle and is thus dissipative (see, e.g., \cite{EthierKurtz1986},
Lemma 4.2.1).

(b)
By Proposition \ref{P3.4}, for any $y\in E$, $x\in E$ and $t>0$, we have
%
%
\begin{equation}
\label{E5.16}
u_y(t)(x)=\mathcal{S}_t F(\bolds{\cdot},y)(x)= F (x,e^{-ct}y
) F \bigl(\theta
, (1-e^{-ct} )y \bigr).
\end{equation}
As $F (\bolds{\cdot},e^{-ct}y )$ is in $C_l^2(E)$, so is
$\mathcal{S}_t F(\bolds{\cdot},y)$.

(c)
By (\ref{E5.10}), we have
\[
{\mathcal{G}}^{c,\theta}u_y(t)(x)=\lim_{{\varepsilon}\downarrow0}
{\varepsilon}^{-1} (\mathcal{S}_{t+{\varepsilon}}-\mathcal
{S}_t)F(\bolds{\cdot},y)(x)
=\frac{d}{dt} (u_y(t)(x) ).
\]
Using (\ref{E5.16}) and Lemma \ref{L2.2}, for every $x\in E$, we get
\begin{eqnarray*}
{\mathcal{G}}^{c,\theta}\mathcal{S}_t F(\bolds{\cdot},y)(x)
&=&\langle-c e^{-ct}y,\nabla_2 F(c,e^{-ct}y)\rangle F\bigl(\theta
,(1-e^{-ct})y\bigr)\\
&&{} +F(x,e^{-ct}y)\bigl\langle c e^{-ct}y,\nabla_2F\bigl(\theta
,(1-e^{-ct})y\bigr)\bigr\rangle\\
&=& [ c e^{-ct} (\theta-x)\diamond y]  F (x,e^{-ct}y
)F \bigl(\theta, (1-e^{-ct} )y \bigr).
\end{eqnarray*}
Hence, $t\mapsto{\mathcal{G}}^{c,\theta}u_y(t)$ is clearly
continuous [in $C_l(E)$].

(d)
As $t\mapsto{\mathcal{G}}^{c,\theta}u_y(t)$ is continuous, it is
integrable, and
\begin{eqnarray*}
\biggl(\int_0^t{\mathcal{G}}^{c,\theta}u_y(s) \,ds \biggr)(x)=\int
_0^t{\mathcal{G}}
^{c,\theta}u_y(s)(x)\, ds=u_y(t)(x)-u_y(0)(x)
\end{eqnarray*}
implies (\ref{E5.15}).

\subsection[The martingale problem (proof of Theorem 5.4)]{The martingale
problem (proof of Theorem \protect\ref{T5.4})}
\label{A5.3}

Before we prove this theorem, we derive a duality relation for
processes satisfying the martingale problem $({\mathcal{G}}^{c,\theta},V)$.
Recall the definition of $\tilde Y$ from (\ref{E2.11}).
\begin{lemma}
\label{L5.7}
Let $\mu$ be a probability measure on $E$.
Let $X$ be any solution of the martingale problem $({\mathcal
{G}}^{c,\theta},V)$
with ${\mathcal{L}}[X_0]=\mu$.
Then $X$ and $\tilde Y$ are dual, in the sense that
for any $\tilde y\in E\times[0,\infty)^2$, we have
\[
\mathbf{E}_{\mu} [H ((X_t,\theta),\tilde y ) ]
=\int_E\mathbf{E}_{\tilde y} [H ((x,\theta),\tilde Y_t
)
]\mu(dx) \qquad\mbox{for all }t\geq0.
\]
\end{lemma}
\begin{pf}
As $X$ is a solution of the martingale problem, we have that
\begin{eqnarray*}
&& H ((X_t,\theta),\tilde y )-\int_0^tH ((X_s,\theta
),\tilde y ) [c(\theta-X_s)\diamond\tilde y(1)] \,ds\\
&&\qquad =F(\theta,\tilde y(2)) \biggl(
F(X_t,\tilde y(1))-\int_0^t F(X_s,\tilde y(1)) [c(\theta-X_s)\diamond
\tilde y(1)] \,ds \biggr)
\end{eqnarray*}
is a martingale.
On the other hand, by (\ref{E2.16}) [since $\tilde y(1)\in E$, one term
vanishes],
%
%
\begin{eqnarray}
\label{E5.17}
\frac{d}{dt}\mathbf{E}_{\tilde y}
 [H ((x,\theta),\tilde Y_t ) ]|_{t=0}
&=&
 \langle-c \tilde y(1), \nabla_2F(x,\tilde y(1)) \rangle
F(\theta,\tilde y(2))\nonumber\\
&&{} + \langle c \tilde y(1), \nabla_2F(\theta,\tilde
y(2)) \rangle F(x,\tilde y(1))\\
&=& H ((\theta,x),\tilde y ) [c(\theta-x)\diamond\tilde y(1)].\nonumber
\end{eqnarray}
Since $\tilde Y$ is deterministic, we get that
\[
H ((x,\theta),\tilde Y_t )-\int_0^tH ((\theta,x),\tilde
Y_s ) c(\theta-x)\diamond\tilde Y_s(1) \,ds = H ((x,\theta
),\tilde y )
\]
is the trivial martingale. By \cite{EthierKurtz1986},
Corollary 4.4.13, this implies
that
\[
\mathbf{E}_\mu [H ((X_t,\theta),\tilde y ) ]
=\int\mathbf{E}_y [H ((x,\theta),\tilde Y_t ) ] \mu(dx)
\]
and we are done.
\end{pf}
\begin{pf*}{Proof of Theorem \ref{T5.4}}
By Theorem \ref{T5.3}(ii) and (\ref{E5.11}), and since $V\subset
C_l(E)$, by definition of $X^{c,\theta}$, the process $X^{c,\theta}$
is, in fact, a solution of the martingale problem $({\mathcal
{G}}^{c,\theta}, V)$.

Now, assume that $X$ and $X'$ are two solutions with ${\mathcal
{L}}[X_0]={\mathcal{L}}
[X'_0]=\mu$.
By Lemma \ref{L5.7}, we
get
\[
\mathbf{E}_{\mu} [F (X_t,y ) ]
=\mathbf{E}_{\mu} [F (X'_t,y ) ] \qquad\mbox{for all
}t\geq0\mbox{ and } y\in E.
\]
By Corollary \ref{C2.4}, $\{F(\bolds{\cdot},y), y\in E\}$
is measure determining on $E$. Hence, ${\mathcal{L}}_\mu
[X_t]={\mathcal{L}}_\mu[X_t']$
for all $t\geq0$.
By \cite{EthierKurtz1986}, Theorem 4.4.2, this implies that the
finite-dimensional distributions of $X$ and $X'$ coincide.
\end{pf*}

\section*{Acknowledgment}
The authors would like to express their gratitude to an anonymous
referee for valuable comments.


%
\printaddresses


\begin{thebibliography}{16}

\bibitem{Burkholder1977}
\begin{barticle}[mr]
\bauthor{\bsnm{Burkholder},~\bfnm{D.~L.}\binits{D.~L.}}
(\byear{1977}).
\btitle{Exit times of {B}rownian motion, harmonic majorization, and {H}ardy
  spaces}.
\bjournal{Adv. Math.}
\bvolume{26}
\bpages{182--205}.
\bid{mr={0474525}}
\end{barticle}
\endbibitem

\bibitem{CoxDawsonGreven2004}
\begin{barticle}[mr]
\bauthor{\bsnm{Cox},~\bfnm{J.~T.}\binits{J.~T.}},
  \bauthor{\bsnm{Dawson},~\bfnm{D.~A.}\binits{D.~A.}} \AND
  \bauthor{\bsnm{Greven},~\bfnm{A.}\binits{A.}}
(\byear{2004}).
\btitle{Mutually catalytic super branching random walks: Large finite systems
  and renormalization analysis}.
\bjournal{Mem. Amer. Math. Soc.}
\bvolume{171}
\bpages{viii$+$97}.
\bid{mr={2074427}}
\end{barticle}
\endbibitem

\bibitem{CoxKlenkePerkins2000}
\begin{bincollection}[mr]
\bauthor{\bsnm{Cox},~\bfnm{J.~Theodore}\binits{J.~T.}},
  \bauthor{\bsnm{Klenke},~\bfnm{Achim}\binits{A.}} \AND
  \bauthor{\bsnm{Perkins},~\bfnm{Edwin~A.}\binits{E.~A.}}
(\byear{2000}).
\btitle{Convergence to equilibrium and linear systems duality}.
In \bbooktitle{Stochastic Models ({O}ttawa, {ON}, 1998)}.
\bseries{CMS Conference Proceedings}
\bvolume{26}
\bpages{41--66}.
\bpublisher{Amer. Math. Soc.}, \baddress{Providence, RI}.
\bid{mr={1765002}}
\end{bincollection}
\endbibitem

\bibitem{DawsonFleischmannXiong2005}
\begin{barticle}[mr]
\bauthor{\bsnm{Dawson},~\bfnm{Donald~A.}\binits{D.~A.}},
  \bauthor{\bsnm{Fleischmann},~\bfnm{Klaus}\binits{K.}} \AND
  \bauthor{\bsnm{Xiong},~\bfnm{Jie}\binits{J.}}
(\byear{2005}).
\btitle{Strong uniqueness for cyclically symbiotic branching diffusions}.
\bjournal{Statist. Probab. Lett.}
\bvolume{73}
\bpages{251--257}.
\bid{doi={10.1016/j.spl.2005.03.012}, mr={2179284}}
\end{barticle}
\endbibitem

\bibitem{DawsonPerkins1998}
\begin{barticle}[mr]
\bauthor{\bsnm{Dawson},~\bfnm{Donald~A.}\binits{D.~A.}} \AND
  \bauthor{\bsnm{Perkins},~\bfnm{Edwin~A.}\binits{E.~A.}}
(\byear{1998}).
\btitle{Long-time behavior and coexistence in a mutually catalytic branching
  model}.
\bjournal{Ann. Probab.}
\bvolume{26}
\bpages{1088--1138}.
\bid{doi={10.1214/aop/1022855746}, mr={1634416}}
\end{barticle}
\endbibitem

\bibitem{EthierKurtz1986}
\begin{bbook}[mr]
\bauthor{\bsnm{Ethier},~\bfnm{Stewart~N.}\binits{S.~N.}} \AND
  \bauthor{\bsnm{Kurtz},~\bfnm{Thomas~G.}\binits{T.~G.}}
(\byear{1986}).
\btitle{Markov Processes: Characterization and Convergence}.
\bpublisher{Wiley}, \baddress{New York}.
\bid{doi={10.1002/9780470316658}, mr={838085}}
\end{bbook}
\endbibitem

\bibitem{IkedaWatanabe1989}
\begin{bbook}[mr]
\bauthor{\bsnm{Ikeda},~\bfnm{Nobuyuki}\binits{N.}} \AND
  \bauthor{\bsnm{Watanabe},~\bfnm{Shinzo}\binits{S.}}
(\byear{1989}).
\btitle{Stochastic Differential Equations and Diffusion Processes},
\bedition{2nd} ed.
\bseries{North-Holland Mathematical Library}
\bvolume{24}.
\bpublisher{North-Holland}, \baddress{Amsterdam}.
\bid{mr={1011252}}
\end{bbook}
\endbibitem

\bibitem{Klenke2008e}
\begin{bbook}[mr]
\bauthor{\bsnm{Klenke},~\bfnm{Achim}\binits{A.}}
(\byear{2008}).
\btitle{Probability Theory}.
\bpublisher{Springer}, \baddress{London}.
\bid{doi={10.1007/978-1-84800-048-3}, mr={2372119}}
\end{bbook}
\endbibitem

\bibitem{KM2}
\begin{bmisc}[auto:SpringerTagBib|2009-01-14|16:51:27]
\bauthor{\bsnm{Klenke},~\bfnm{Achim}\binits{A.}} \AND
\bauthor{\bsnm{Mytnik},~\bfnm{Leonid}\binits{L.}}
(\byear{2008}).
\bhowpublished{Infinite rate mutually catalytic
  branching in infinitely many colonies. {C}onstruction, characterization and
  convergence. Preprint. Available at}
\href{http://arxiv.org/abs/arXiv:0901.0623}{arXiv:0901.0623 [math.PR]}.
\end{bmisc}
\endbibitem

\bibitem{KM3}
\begin{bmisc}[auto:SpringerTagBib|2009-01-14|16:51:27]
\bauthor{\bsnm{Klenke},~\bfnm{Achim}\binits{A.}} \AND
\bauthor{\bsnm{Mytnik},~\bfnm{Leonid}\binits{L.}}
(\byear{2009}).
\bhowpublished{Infinite rate mutually catalytic
  branching in infinitely many colonies. {T}he longtime behaviour.
  Preprint. Available at}
\href{http://arxiv.org/abs/arXiv:0910.4120}{arXiv:0910.4120 [math.PR]}.
\end{bmisc}
\endbibitem

\bibitem{KO}
\begin{barticle}[vtex]
\bauthor{\bsnm{Klenke},~\bfnm{Achim}\binits{A.}} \AND
\bauthor{\bsnm{Oeler},~\bfnm{Mario}\binits{M.}}
(\byear{2010}).
\btitle{A Trotter type approach to
  infinite rate mutually catalytic branching}.
\bjournal{Ann. Probab}.
\bvolume{38}
\bpages{479--497}.
\end{barticle}
\endbibitem

\bibitem{LeGallMeyre1992}
\begin{barticle}[mr]
\bauthor{\bsnm{Le~Gall},~\bfnm{J.~F.}\binits{J.~F.}} \AND
  \bauthor{\bsnm{Meyre},~\bfnm{T.}\binits{T.}}
(\byear{1992}).
\btitle{Points c\^ones du mouvement brownien plan, le cas critique}.
\bjournal{Probab. Theory Related Fields}
\bvolume{93}
\bpages{231--247}.
\bid{doi={10.1007/BF01195230}, mr={1176726}}
\end{barticle}
\endbibitem

\bibitem{Mytnik1998}
\begin{barticle}[mr]
\bauthor{\bsnm{Mytnik},~\bfnm{Leonid}\binits{L.}}
(\byear{1998}).
\btitle{Uniqueness for a mutually catalytic branching model}.
\bjournal{Probab. Theory Related Fields}
\bvolume{112}
\bpages{245--253}.
\bid{doi={10.1007/s004400050189}, mr={1653845}}
\end{barticle}
\endbibitem

\bibitem{Oeler2008}
\begin{bmisc}[auto:SpringerTagBib|2009-01-14|16:51:27]
\bauthor{\bsnm{Oeler},~\bfnm{Mario}\binits{M.}}
(\byear{2008}).
\bhowpublished{Mutually catalytic branching at infinite
rate. Ph.D. thesis, Univ. Mainz}.
\end{bmisc}
\endbibitem

\bibitem{RogersWilliams1994}
\begin{bbook}[mr]
\bauthor{\bsnm{Rogers},~\bfnm{L.~C.~G.}\binits{L.~C.~G.}} \AND
  \bauthor{\bsnm{Williams},~\bfnm{David}\binits{D.}}
(\byear{1994}).
\btitle{Diffusions, {M}arkov Processes, and Martingales. {V}ol. 1},
\bedition{2nd} ed.
\bpublisher{Wiley}, \baddress{Chichester}.
\bid{mr={1331599}}
\end{bbook}
\endbibitem

\bibitem{ShigaShimizu1980}
\begin{barticle}[mr]
\bauthor{\bsnm{Shiga},~\bfnm{Tokuzo}\binits{T.}} \AND
  \bauthor{\bsnm{Shimizu},~\bfnm{Akinobu}\binits{A.}}
(\byear{1980}).
\btitle{Infinite-dimensional stochastic differential equations and their
  applications}.
\bjournal{J. Math. Kyoto Univ.}
\bvolume{20}
\bpages{395--416}.
\bid{mr={591802}}
\end{barticle}
\endbibitem

\end{thebibliography}
\end{document}